\documentclass[11pt]{amsart}

\usepackage{amssymb}
\usepackage{a4wide}

\newtheorem{theorem}{Theorem}
\newtheorem{prop}{Proposition}
\newtheorem{lemma}{Lemma}
\newtheorem{coro}{Corollary}
\newtheorem{fact}{Fact}
 
\newcommand{\ts}{\hspace{0.5pt}}
\newcommand{\CC}{\mathbb{C}\ts}
\newcommand{\RR}{\mathbb{R}\ts}
\newcommand{\QQ}{\mathbb{Q}\ts}

\newcommand{\ZZ}{\mathbb{Z}}
\newcommand{\NN}{\mathbb{N}}

\newcommand{\cG}{\mathcal{G}}
\newcommand{\cA}{\mathcal{A}}
\newcommand{\cE}{\mathcal{E}}
\newcommand{\cB}{\mathcal{B}}
\newcommand{\cT}{\mathcal{T}}
\newcommand{\cI}{\mathcal{I}}
\newcommand{\cJ}{\mathcal{J}}
\newcommand{\cH}{\mathcal{H}}
\newcommand{\cF}{\mathcal{F}}
\newcommand{\cS}{\mathcal{S}}
\newcommand{\cR}{\mathcal{R}}
\newcommand{\dd}{{\rm d}}

\DeclareMathOperator{\len}{len}
\DeclareMathOperator{\lcm}{lcm}
\DeclareMathOperator{\ord}{ord}
\DeclareMathOperator{\trace}{tr}
\DeclareMathOperator{\cent}{cent}
\DeclareMathOperator{\ch}{char}
\DeclareMathOperator{\diag}{diag}
\DeclareMathOperator{\pd}{pol\,deg}
\newcommand{\amalgam}{\hspace{1pt}\begin{smallmatrix} \ast \\ 
      \!{\scriptscriptstyle \cB} \end{smallmatrix}}
\newcommand{\bs}[1]{\boldsymbol{#1}}

\begin{document}

\title[Symmetries and reversing symmetries]
{Symmetries and reversing symmetries of \\[1mm]
polynomial automorphisms of the plane}

\author{Michael Baake}
\address{Fakult\"at f\"ur Mathematik, Universit\"at Bielefeld, 
Box 100131, 33501 Bielefeld, Germany}
\email{mbaake@math.uni-bielefeld.de}

\author{John A.~G.~Roberts}
\address{School of Mathematics, 
University of New South Wales,
Sydney, NSW 2052, Australia}
\email{jag.roberts@unsw.edu.au}

\begin{abstract} 
The polynomial automorphisms of the affine plane over a field $K$ form
a group which has the structure of an amalgamated free product.  This
well-known algebraic structure can be used to determine some key
results about the symmetry and reversing symmetry groups of a given
polynomial automorphism.
\end{abstract}

\maketitle

\section{Introduction}

In a series of recent articles \cite{BRcat,BRtorus,RBtrace,rowi}, the
symmetries and reversing symmetries of some dynamical systems
(automorphisms) have been investigated systematically by means of
algebraic methods.  An automorphism $L$ of some space is said to have
a {\em symmetry\/} if there exists an automorphism $S$ that satisfies
\begin{equation} \label{defsym}
        S \circ L \circ S^{-1} \; = \; L \ts ,
\end{equation}
and a {\em reversing symmetry\/}, or {\em reversor}, if there exists an 
automorphism $R$ so that
\begin{equation} \label{defrevsym}
         R \circ L \circ R^{-1} \; = \;  L^{-1}\ts .
\end{equation}
The set of symmetries is non-empty (it certainly contains all powers
of $L$) and this set is actually a group, the {\em symmetry group\/}
$\cS(L)$. On the other hand, the existence {\em a priori\/} of any
reversing symmetries for a particular $L$ is unclear. When $L$ has a
reversing symmetry, we call it {\em reversible}, and {\em
irreversible} otherwise.  The set $\cR(L)$ of all symmetries and
reversing symmetries of $L$ is a group, too, called the {\em reversing
symmetry group\/} \cite{lamb} of $L$ (see also \cite{Goodson}).  In
particular, $\cR(L)$ admits a binary grading: the composition of two
reversing symmetries is a symmetry, whereas the composition of a
symmetry and a reversing symmetry is a reversing symmetry. If $L$ is
irreversible or if $L$ is the identity or an involution (i.e., if $L^2
= 1$), one has $\cR(L)=\cS(L)$; otherwise, $\cR(L)$ is a group
extension of $\cS(L)$ of index $2$.

The simultaneous consideration of symmetries and reversing symmetries
of reversible automorphisms (which may arise as the time-one maps of
reversible flows) is now known to provide some powerful algebraic
insights. For example, the results of \cite{lamb,Goodson} illustrate
that much can be said about the nature of possible reversing
symmetries in $\cR(L)$ given the knowledge of the structure of
$\cS(L)$. For example, if $L$ (with $L^2\neq 1$) has an involutory
reversor $R$ (i.e., $R^2 = 1 \neq R$), one has $\cR(L)\simeq \cS(L)
\rtimes C_2$,
where $C_n$ denotes the cyclic group of order $n$ and
$N \rtimes H$ is the semi-direct product of $N$ and $H$,
with $N$ the normal subgroup.
In many cases of reversible automorphisms (and also in the analogous
continuous-time case of reversible flows), it is in fact found that
all reversing symmetries $R$ that satisfy (\ref{defrevsym}) are
involutions. In this case, the automorphism $L$ can be written as the
composition of two involutions, e.g., $L \circ R$ and $R$, or $R$ and
$R \circ L$. References \cite{roqu} and \cite{LR} include reviews of
the properties and applications of reversible automorphisms and flows.

The programme followed in the papers \cite{BRcat,BRtorus,RBtrace,rowi}
can be summarised as follows.  The nature of $\cS(L)$ and $\cR(L)$ has
been investigated for some well-known groups of automorphisms where
the group structure admits an algebraic investigation of the relations
(\ref{defsym}) and (\ref{defrevsym}). This necessitates restricting
the search for $S$ and $R$ to some suitable group that contains the
automorphism $L$ (which might be argued to be a natural first step).
Dynamical systems considered in this programme have included toral
automorphisms in two and higher dimensions and polynomial automorphism
of $\RR^3$ that are closely related (by semi-conjugacies) to
two-dimensional toral automorphisms and arise as {\em trace maps\/} in
the study of quasi-periodic phenomena and the theory of aperiodic
order.

In \cite{RBstandard}, we turned our attention to the group of planar
polynomial automorphisms. This group comprises ``maps'' of the form 
\begin{equation} \label{defmap} 
     x' \; = \; P(x,y) \; , \quad y' \; = \; Q(x,y)\, ,  
\end{equation}
where $P(x,y)$ and $Q(x,y)$ are polynomials with coefficients in some
field $K$, {\em and} there is an inverse that is also polynomial (so
the polynomial map $x' = x^3, \, y'=x+y$, although a bijection over
$\RR^2$, is {\em not} in the group since its inverse involves cube
roots; see \cite{Rudin} for the contrasting complex case). The term
``map'' in this context is actually a slight abuse of language. Over
finite fields, different polynomials (such as $P_1(x)=x$ and
$P_2(x)=x^2$ over $\mathbb{F}_2$, the finite field with two elements)
can define the same mapping. If we use the term ``map'' or
``polynomial map'' in this article, we actually mean to distinguish
them according to their polynomial structure.

The group of polynomial automorphisms of the plane $K^2$, denoted
${\rm GA}^{}_2(K)$, has been studied in some detail because it has the
structure of an amalgamated free product (compare \cite{Essen,Essen2}
and references therein, and Section 2 below for more
details). Obviously, polynomial maps are much-studied as dynamical
systems. In particular, ${\rm GA}^{}_2(\RR)$ and ${\rm GA}^{}_2(\CC)$
have received considerable attention. They include, for example, the
H\'enon quadratic map family,
\begin{equation} \label{defhenonmap} 
     x' \; = \; y  \; , \quad y' \; = \; -\delta\ts x + y^2 + c\ts ,  
\end{equation}
with constants $c,\delta\in\CC$ and $\delta\neq 0$. This is one of the
more famous ``toy models'' of discrete dynamics.  Exploitation of the
group structure of ${\rm GA}^{}_2(\RR)$ and ${\rm GA}^{}_2(\CC)$ has
been used to great effect to investigate various properties of their
elements, e.g., their roots \cite{AR} or their dynamical entropy
\cite{FM}.  The same idea was used in \cite{RBstandard} to give a
description of possible $\cS(L)$ and $\cR(L)$ structures for
the subset of ${\rm GA}^{}_2(\RR)$ of maps in so-called {\em
generalised standard form}
\begin{equation} \label{defgenstan} 
     x' \; = \; x+P_1(y)  \; , \quad y' \; = \; y + P_2(x') \, ,  
\end{equation}
with polynomials $P_1$ and $P_2$ (and inverse: $y = y' - P_2(x')$, 
$x = x' - P_1(y)$).  The form (\ref{defgenstan}) is a common one for
area-preserving maps in the dynamics literature. In \cite{RBstandard},
we also provided normal forms for maps of the form (\ref{defgenstan})
with the various possible symmetries or reversing symmetries.
Subsequently, G\'omez and Meiss \cite{GM} have given normal forms for
general elements of ${\rm GA}^{}_2(\RR)$ and ${\rm GA}^{}_2(\CC)$ that
possess involutory reversing symmetries.

In this paper, we return to symmetries and reversing symmetries of
general elements of ${\rm GA}_2(K)$. As compared to \cite{RBstandard}
and \cite{GM}, our approach will be significantly more algebraic.
This is possible due to ${\rm GA}_2(K)$ being an amalgamated free
product of two well understood groups, so that combinatorial group theory
can be used very effectively. Unfortunately, since no such structure
is at hand for more than two dimensions, compare \cite[Ex.\
2.4]{Kraft}, our approach is presently restricted to the planar case.

Particular goals of this paper are:\ (i) to make maximal use of the
algebraic consequences of the amalgamated free product structure of
the group; (ii) to concentrate on characterising $\cS(L)$ {\em before}
moving onto the study of reversing symmetries, in view of the benefits
that can flow algebraically in this direction; and (iii) to carry
through some of the results for a general field $K$, before
specialising to $\RR$ or $\CC$. Of course, the real and complex cases
would seem to be the most interesting ones historically. However,
dynamical systems over finite fields are becoming more topical, see
\cite{rovi1,rovi2} and references therein. In particular, \cite{rovi2}
studies the cycle statistics of permutations associated with
reductions to finite fields of planar polynomial automorphisms. It
turns out that application to the finite fields case of
Proposition~\ref{elementary-to-linear} of Section 4 below helps to
understand how the possession of orientation-reversing involutory
polynomial reversors leads to more cycles of shorter average length than would
otherwise occur in, for example, random permutations.

As an indication of the results we obtain via this algebraic approach,
we mention some of them for a ``typical'' infinite order element $L$
of ${\rm GA}^{}_2(\RR)$ (here, ``typical'' means a {\em CR element}
$f$, cf.\ Section 2):
\begin{itemize}
\item any nontrivial symmetry of $L$ of finite order is an involution
conjugate to $\left(\begin{smallmatrix} -1 & 0\\ 0 & -1
\end{smallmatrix}\right)$
and $\cS(L)$ contains at most one nontrivial finite subgroup (then
isomorphic with the cyclic group $C_2$) [Theorem \ref{involutions},
Section \ref{symmetries}]; a strong characterisation can be given both
for the involutory symmetry and for $L$ [Theorem \ref{symnormform} and
Corollary \ref{involtest}, Section \ref{symmetries}];
\item any reversor of $L$ is of finite order,
being an involution, conjugate to
$\left(\begin{smallmatrix} -1 & 0\\ 0 & -1
\end{smallmatrix}\right)$
or
$\left(\begin{smallmatrix} 0 & 1\\ 1 & 0
\end{smallmatrix}\right)$,
or an element of order 4, conjugate to 
$\left(\begin{smallmatrix} 0 & -1 \\ 1 & 0
\end{smallmatrix}\right)$
[Theorems \ref{rev-order-real} and \ref{natureorder4}, 
Section \ref{revsymmsec}];
\item if $L$ has a reversor, a normal form to which $L$ is conjugate
in ${\rm GA}^{}_2(\RR)$ can be found [Propositions \ref{invnormform}
and \ref{4normform}, Section \ref{revsymmsec}].
\end{itemize}

The plan of the paper is as follows.  In Section 2, we recall key
results about the group structure of the planar polynomial
automorphisms.  In Section 3, we summarise results from \cite{Wright}
and various other sources \cite{MKS,LS,Mol} concerning Abelian
subgroups of ${\rm GA}^{}_2(K)$.  This is exploited in Section 4,
where we characterise the symmetry groups and symmetries of typical
elements.  Finally, in Section 5, we employ knowledge of the
symmetries to characterise the possible reversing symmetries.

\bigskip
In the final preparation of this manuscript, we became aware of
related results by Gom\'ez and Meiss in a preprint that has now
appeared \cite{GM2}. They concentrate on the cases $K=\RR$ and
$K=\CC$, using rather explicit calculations with normal forms, while
our focus is more on the general setting, with stronger focus on
algebraic methods.  We compare their main results with ours in remarks
preceding Theorem 2 in Section 4 and following Propositions 11 and 12
in Section 5, augmented by various smaller remarks throughout the
paper.

\section{Recollections and mathematical setting}

Let us first recall a number of well-known results about the group
structure of the polynomial automorphisms of the plane.  We do this in
some generality, and simultaneously introduce our notation. Most of
what is contained in this section is classic material, and mainly
relies on \cite{Wright,FM,KS} and references given therein. Still, it
seems worthwhile to combine several results in a fashion that suits
our purpose and makes the paper more self-contained.

\smallskip
Let $K$ be a field and consider the group $\cG^{}_K = {\rm
GA}^{}_2(K)$ of polynomial automorphisms of the affine plane over $K$,
i.e., the set of mappings of the form \eqref{defmap} with $P,Q\in
K[x,y]$ (the ring of polynomials in $x,y$ with coefficients in $K$)
such that the inverse exists and is also polynomial.  Group
multiplication is composition of maps, usually written as $gg'$ rather
than $g\circ g'$ in the sequel. The neutral element of the group will
be written as $1$, denoting the identity map. For mappings, in
comparison with \eqref{defmap}, we interchangeably also use the
notation
\begin{equation} \label{polydef}
    \binom{x}{y} \; \mapsto \; \binom{P(x,y)}{Q(x,y)} .
\end{equation}
Whenever $K$ is clear from the context, we will write $\cG$ rather
than $\cG^{}_K$ for simplicity.

Note that the Jacobian ${\rm d}g$ of any element $g\in\cG$ (defined
via the algebraic derivative of the polynomials involved, see
\cite[p.~5]{Essen}) has constant determinant $\neq 0$, i.e.,
$\det({\rm d}g)$ is an element of $K^* =K\setminus \{0\}$, the latter
representing the only units of the ring $K[x,y]$, compare
\cite[Prop.~1.14]{Essen}. The converse question is connected with the
famous Jacobian conjecture, namely whether $\det({\rm d}g)=const\neq
0$ is sufficient for a polynomial mapping to be an automorphism, see
\cite{Essen} for a summary and \cite{Rudin} for an interesting partial
result, together with some comments on the influence of the field $K$
being algebraically closed or not.

The group $\cG$ contains two particularly important subgroups. First,
there is the group $\cA$ of {\em affine}\/ transformations,
\begin{equation*} 
     \cA \; = \; \{ (\bs{a},M) \mid \bs{a}\in K^2, \, M\in{\rm GL}(2,K)\}, 
\end{equation*}
where $(\bs{a},M)$ encodes the mapping $\bs{x}\mapsto M\bs{x} +
\bs{a}$. We write $\bs{x}$ for a column vector with two entries, and
tacitly identify the elements of $\cA$ with the corresponding ones of
$\cG$.  In particular, a matrix $M$ is identified with the linear
mapping $\bs{x}\mapsto M\bs{x}$, and a vector space element $\bs{a}\in
K^2$ with the translation $\bs{x}\mapsto\bs{x}+\bs{a}$.  The
multiplication of two elements of $\cA$ reads
\begin{equation*} 
    (\bs{a},A) (\bs{b},B) \; = \; (\bs{a} + A \bs{b}, AB) 
\end{equation*}
which shows that $\cA$ is a semi-direct product, $\cA = K^2 \rtimes
{\rm GL}(2,K)$, where $K^2$ is the normal subgroup. Note that the
inverse of $(\bs{a},A)$ reads $(\bs{a},A)^{-1} = (-A^{-1}\bs{a},
A^{-1})$.

The second group, $\cE$, consists of all mappings of $\cG$ of the form
\begin{equation} \label{elementary}
    e : \quad
    \binom{x}{y} \; \mapsto \; \binom{\alpha x + P(y)}{\beta y + v} 
\end{equation}
with $P$ a polynomial, $\alpha,\beta,v \in K$ and $\alpha \beta \neq 0$. 
The inverse reads
\begin{equation} \label{inv-elementary}
   e^{-1} : \quad
   \binom{x}{y} \; \mapsto \; 
   \binom{\frac{1}{\alpha}\, (x - P(\frac{y-v}{\beta})) }
         {\frac{1}{\beta}\, (y - v)} . 
\end{equation}
The elements of $\cE$ are called {\em elementary}\/ transformations.
They map lines with constant $y$-coordinate to lines of the same type.
The relevance of these two subgroups comes from the following fact,
which was proved by Jung \cite{Jung} for $K\in\{\RR,\CC\}$ and later
by van der Kulk \cite{Kulk} for arbitrary fields $K$, see also
\cite[Sec.~1.5]{Wright} and \cite[p.~68]{FM}.
\begin{fact} \label{generating}
  The group\/ $\cG$ of polynomial automorphisms of the plane\/ $K^2$
  is generated by the two subgroups\/ $\cA$ and\/ $\cE$.  \qed
\end{fact}

The intersection of $\cA$ and $\cE$, both seen as subgroups of $\cG$,
is another group, called $\cB$ (for {\em basic}) from now on. It
consists of all mappings of the form
\begin{equation} \label{basic}
   b \; : \quad \binom{x}{y} \; \mapsto \; 
   \begin{pmatrix} \alpha & \gamma \\ 0 & \beta \end{pmatrix}
   \binom{x}{y} + \binom{u}{v}    
\end{equation}
with $\alpha,\beta,\gamma,u,v\in K$ and $\alpha \beta \neq 0$.  If
$\cT$ denotes the subgroup of ${\rm GL}(2,K)$ which consists of all
upper (invertible) triangular matrices, one can see that $\cB$ is
again a semi-direct product,
\begin{equation*} 
    \cB \; = \; K^2 \rtimes \cT \, .
\end{equation*}
The following result \cite{Serre,Wright} is important, see
also \cite[Thm.~5.1.11 and Cor.~5.3.6]{Essen2}.
\begin{fact} \label{amalgamation}
   The group\/ $\cG$ is the free product of the groups\/ $\cA$ and\/
   $\cE$, called factors, amalgamated along their intersection, $\cB$,
   abbreviated as\/ $\cG = \cA \amalgam \cE$.  \qed
\end{fact}

This gives access to the structure of the group $\cG$, and to its
subgroups in particular. To make explicit use of it later on, we need
a natural way to represent group elements uniquely. This is achieved
by a partition of $\cG$ into (right) cosets, which we will indicate
below by the symbol $\dot{\cup}$ (for disjoint union). Define
\begin{equation} \label{I-def}
   \cI \; := \; \left\{ \begin{pmatrix} 0 & 1 \\ 1 & \beta \end{pmatrix}
   \Big\vert\; \beta \in K \right\}  \; \subset \; \cA\setminus\cB \, , 
\end{equation}
again identified with the corresponding subset of $\cG$, and
\begin{equation} \label{J-def}
   \cJ \; := \; \left\{ \binom{x}{y}\mapsto
   \binom{x+y^2P(y)}{y}\, \Big\vert\; 0\neq P\in K[y]  \right\} 
   \; \subset \; \cE\setminus\cB \, . 
\end{equation}
Note that $\cJ$ is invariant under taking inverses, i.e., if
$e\in\cJ$, so is $e^{-1}$. Furthermore, all elements of $\cI$ and
$\cJ$ fix the origin.

Now, either following \cite[Secs.~1.6 and 1.7]{Wright} and observing
that we use upper triangular matrices in $\cB$ for consistency with
\cite{FM}, or verifying it by a direct computation, one obtains
\begin{fact} \label{cosets}
   Let\/ $\cI$ and\/ $\cJ$ be the sets defined in\/ $(\ref{I-def})$ and\/ 
   $(\ref{J-def})$. Then, the subgroups\/ $\cA$ and\/ $\cE$ of\/ $\cG$ 
   satisfy the unique right coset decompositions\/ 
   $ \cA = \dot{\bigcup}_{a\in\cI\cup\{1\}}\, \cB a $ and\/
   $ \cE = \dot{\bigcup}_{e\in\cJ\cup\{1\}}\, \cB e $
   with respect to the subgroup\/ $\cB=\cA\cap\cE$ of\/ $\cG$.   \qed
\end{fact}

This admits the introduction of a powerful concept, the so-called {\em
normal form}\/ of an element $g\in\cG$, compare \cite[Ch.\
I.1.2]{Serre} and \cite[Sec.\ 4.2]{MKS}.  We also recall the slightly
weaker, but sometimes more useful, result on the {\em reduced word\/}
representation \cite{Cohen,MKS}.
\begin{prop} \label{normal-form}
   Each element\/ $g\in\cG\setminus\cB$ can be written as a\/ {\em reduced
   word}
\begin{equation} \label{reduced-word}
   g \; = \; g^{}_{n} \circ g^{}_{n-1} \circ\ldots\circ g^{}_{1}
\end{equation}
   where\/ $n\ge 1$ and the\/ $g^{}_i$ alternate between\/
   $\cA\setminus\cB$ and\/ $\cE\setminus\cB$, starting and ending with
   either.  Such a product can never be in\/ $\cB$, and it cannot be
   in $\cA\cup\cE$ whenever\/ $n>1$.

   Moreover, each element\/ $g\in\cG$ has a\/ {\em unique}
   representation in the form
\begin{equation} \label{product-form}
   g \; = \; b\circ a^{}_m \circ e^{}_m \circ\ldots\circ a^{}_1 \circ e^{}_1
\end{equation}
   for some\/ $($unique\hspace{1pt}$)$ $m\ge 1$.  Here, $b\in\cB$
   $($including the case that\/ $b$ is the neutral element\/ $1$,
   hence effectively missing\hspace{1pt}$)$, while\/ $a^{}_i \in \cI$
   and\/ $e^{}_i \in \cJ$ for all\/ $1\le i\le m$, except that\/
   $a^{}_m$ and/or\/ $e^{}_1$ are allowed to be missing. The
   representation\/ $\eqref{product-form}$ is called the\/ {\em normal
   form} of\/ $g$ with respect to the coset representatives\/ $\cI$
   and\/ $\cJ$.
\end{prop}
\begin{proof}
The reduced word is a standard way to represent elements in
amalgamated free products, see \cite[Thm.~26]{Cohen}.

The unique normal form emerges as soon as coset representatives of the
factors ($\cA$ and $\cE$ in our case) mod the amalgamation group
($\cB$) are selected. This is achieved by Fact~\ref{cosets}. For
details, see \cite[Thm.~25]{Cohen}.  
\end{proof}

\noindent {\sc Remark}: All elements of $\cI$ (resp.\ $\cJ$)
have Jacobians of determinant $-1$ (resp.\ $+1$), so that
$\lvert\det({\rm d}g)\rvert = \lvert\det({\rm d}b)\rvert$, if
$b\in\cB$ is the starting element according to the decomposition
(\ref{product-form}).  Also, an element in normal form
\eqref{product-form} fixes the origin if and only if the element $b$
in it does.

\medskip
For $g\in\cG$, both the reduced word \eqref{reduced-word} and the
normal form \eqref{product-form}, which we often prefer to deal with,
admit the introduction of several useful concepts, one being the {\em
length}\/ of an element $g$, written as $\len(g)$. If $g$ is given in
normal form (\ref{product-form}), $\len(g)$ is the total number of
factors from $\cI$ and $\cJ$, hence an integer between $2m-2$ and
$2m$, depending on which factors are absent. One has $\len(1)=0$, and,
more generally, $\len(b)=0$ for all $b\in\cB$.
Moreover, the value of $\len(g)$ does not depend on the choice of
coset representatives of the factors, such as $\cI$ and $\cJ$
above. Consequently, if $g\not\in\cB$, $\len(g)$ is the same for all
reduced word representations \eqref{reduced-word} obtained from
\eqref{product-form} by inserting $b^{}_i b^{-1}_i$, with arbitrary
$b^{}_i\in\cB$, between any $a^{}_i$ and the following $e^{}_i$, and
rewriting it in the form \eqref{reduced-word}. One thus has
$\len(g)=n$ for \eqref{reduced-word}.

An element $g\in\cG$ in normal form (\ref{product-form}) is called
{\em cyclically reduced}\/ if $b^{-1} g$ starts with an element in
$\cI$ and ends with one in $\cJ$, or vice versa (so that the word
$b^{-1} g$ alternates between elements of $\cI$ and $\cJ$ when wrapped
on a circle). In other words, $g$ has a cyclically reduced normal form
(CRNF) iff the length of $g$ is even and $>0$. With this definition,
which follows \cite{Serre,Wright,FM} but deviates from \cite{MKS},
elements from the factors or their conjugates cannot have a CRNF, nor
be conjugate to one. This will prove useful shortly.

{}From now on, we will call $g\in\cG$ a {\em CR element}\/ if it is
conjugate to an element with a CRNF. Note that, in general, a CR
element itself does not have a CRNF. As an example, consider the CR
element $g=e^{-1}_1 a^{-1}_1 a e a^{}_1 e^{}_1$, with $a^{-1}_1 a
\not\in\cB$. Then, $g$ is essentially in normal form (possibly with
$b\ne 1$, after rewriting it with the proper representatives from
$\cI$ and $\cJ$), but not cyclically reduced.  Rather, $g$ is
conjugate to $a e$, which is cyclically reduced.

Note that CR elements of $\cG$ play a similar role as hyperbolic
elements do in the class of toral automorphisms \cite{BRtorus}, and
they are the ones we are mainly interested in dynamically. They are
also the ones that can be accessed algebraically, due to the very
structure of $\cG$.  The following result is standard, compare
\cite[Sec.~I.1.3]{Serre} and \cite[Thm.~4.6]{MKS}.
\begin{fact} \label{conjugate}
   Any element\/ $g\in\cG$ is either conjugate to an element of\/
   $\cA$ or\/ $\cE$, or is a CR element, the two cases being mutually
   exclusive.  Moreover, no CR element is of finite order, wherefore
   any element of finite order is conjugate to an element in one of
   the factors.  \qed
\end{fact}
\noindent {\sc Remark}: If $g$ has a CRNF,
one finds $\len(g)\ge 2$ and can check that $\len(g^n) = \lvert n
\rvert \, \len(g)$, for all $n\in\ZZ$ (note that $\len(g^{-1}) =
\len(g)$). Clearly, an element $g$ cannot be of finite order unless
the set $\{\len(g^m)\mid m\ge 0\}$ is bounded and contains $0$. Note,
however, that this is not sufficient for $g$ to be of finite order.
For example, the sequence of lengths is identically $0$ for iterates
of $g\!: x' = c\hspace{0.5pt} x \, , \; y' = c\hspace{0.5pt} y$ where
$c\in K^*$ is an element of infinite multiplicative order (such as
$c=2$ for a field of characteristic $0$).  If all elements of $K^*$
are of finite multiplicative order (e.g., if $K$ is a finite field),
all elements of $\cA$ and $\cE$, and hence all conjugates of such
elements, are of finite order \cite{UH}. In this case, since all
remaining elements are CR elements, finite order and bounded length
are equivalent.

\medskip
{}For completeness, let us recall the following result of Serre
\cite[Sec.\ I.4.3, Thm.\ 8 and its Corollary]{Serre}, which we
formulate in our setting (though it is valid for any amalgamated free
product of two groups).
\begin{fact} \label{serre-fact}
  Any subgroup of\/ $\cG=\cA\amalgam\cE$ of bounded length is
  conjugate to a subgroup of one of the factors.  In particular, every
  finite subgroup of\/ $\cG$ is conjugate to a subgroup of\/ $\cA$
  or\/ $\cE$.
\end{fact}

Note that stronger statements are possible for our special group,
concerning the conjugacy of finite order elements and finite subgroups
to linear ones, if $K$ has characteristic $0$ and is algebraically
closed, compare \cite[p.\ 57 and Thm.\ 2.3]{Kraft} and \cite[Thm.\ 4.3
and Cor.\ 4.4]{Kamba}.  However, this is not so for general $K$, see
\cite{Asa}, wherefore we omit further details.

\medskip
A related concept is that of the {\em degree}\/ of a (non-zero)
polynomial mapping. The degree of $P\in K[x,y]$ is the maximum of the
degrees of its monomials with nonzero coefficient, where $\deg(x^m
y^n)=m+n$, and the degree of the polynomial mapping (\ref{polydef}) is
then defined as the maximum of the degrees of $P$ and $Q$. All affine
maps have degree $1$, and the degree of an elementary map
(\ref{elementary}) is $\max(1,\deg(P))$.  Consequently, the degree
cannot be multiplicative in general, but it is for the decomposition
(\ref{product-form}), see \cite[Thm.~2.1]{FM}.
\begin{fact} \label{multiplicativity}
   If\/ $g\in\cG$ is decomposed according to\/ $\eqref{product-form}$ 
   of Proposition~$\ref{normal-form}$,
   the degree of\/ $g$ is the product of the degrees of the factors, i.e.,
\begin{equation*} 
   \deg(g) \; = \; \prod_{i=1}^{m}\, \deg(e_i) \, , 
\end{equation*}
   where we set\/ $\deg(e^{}_1)=1$ if\/ $e^{}_1$ is missing in the 
   product\/ $\eqref{product-form}$.    \qed
\end{fact}

By analogy to before, the degree of a subgroup is defined as the
maximum of the degrees of its elements. For any $g\in\cG$, one has the
relation $\deg(g) \ge 2^{\lfloor\len(g)/2\rfloor}$, in obvious
modification of \cite[Eq.~(21)]{Wright}. Consequently, the degree
implies a bound on the length, both for elements and for groups. Note,
however, that a group of bounded length need not be of bounded degree;
this can be seen from $\cE$, which is of length $1$, but of unbounded
degree.

\section{Conjugacy and Abelian subgroups in $\cG$}

Let us first recall the following result about conjugacy, where we
rephrase \cite[Thm.~4.6]{MKS} in our terminology. To simplify the
following formulations, we specify:
\begin{itemize}
\item unless stated otherwise, conjugate always means conjugate in $\cG$.
\end{itemize}
Moreover, although ultimately we have our special group $\cG=\cG^{}_K$
in mind, the statements until Theorem~\ref{abelian-subgroups} are not
restricted to this case (unless stated so explicitly), but are
actually valid for the free products of two groups $\cA$ and $\cE$
with an amalgamated subgroup $\cB$.
\begin{prop} \label{conjugacy}
   In the amalgamated free product\/ $\cG=\cA\amalgam\cE$,
   every element of\/ $\cG$ is conjugate to an element of\/ $\cA$ or\/ 
   $\cE$, or is a CR element, i.e., conjugate to an element with CRNF.

   Moreover, if\/ $g$ is itself an element of\/ $\cA\cup\cE$ or 
   an element with CRNF, one has the following three possibilities.
\begin{enumerate}
\item If\/ $g$ is conjugate to an element\/ $b\in\cB$, then\/ $g$ lies in
   one of the factors, and there is a sequence\/ $b, h^{}_1, h^{}_2,
   \dots, h^{}_{\ell}, g$ where each\/ $h_i$ lies in\/ $\cB$ and
   consecutive elements of the sequence are conjugate in a factor.
\item If\/ $g$ is conjugate to an element\/ $g'$ that is in some factor 
   but not in a conjugate of\/ $\cB$, then also\/ $g$ lies in the same 
   factor, and\/ $g$ and\/ $g'$ are conjugate within this factor.
\item If\/ $g$ is conjugate to an element\/ $g'$ in CRNF,
   one can obtain\/ $g$ from\/ $g'$ by a cyclic 
   permutation of the factors of\/ $g'$, followed by a conjugacy
   with an element from\/ $\cB$.   
\end{enumerate}
For general elements, these possibilities apply up to conjugacy. \qed
\end{prop}

An amalgamated free product of two factors admits some access to the 
structure of its subgroups, in particular the Abelian ones. Let us 
first consider two commuting elements.
\begin{lemma} \label{commute}
   Let\/ $g,g'\in\cG$ with\/ $g\hspace{1pt}g' = g'g$. Then, one of the 
   following three cases applies.
\begin{enumerate}
\item The element\/ $g$ or\/ $g'$ is in a conjugate of\/ $\cB$.
\item If neither\/ $g$ nor\/ $g'$ is in a conjugate of\/ $\cB$, but\/
      $g$ is in a conjugate of a factor, then\/ $g'$ is in that
      same conjugate, too.
\item If neither\/ $g$ nor\/ $g'$ is in a conjugate of a factor\/
      $($i.e., if both\/ $g$ and\/ $g'$ are CR elements\/$)$, then\/
      $g =d^k c\ts b c^{-1}$ and\/ $g' = d^{\ell} c\ts b' c^{-1}$ for
      some\/ $c,d\in\cG$, $k,\ell\in\ZZ$ and\/ $b,b'\in\cB$, where\/
      $c\ts b c^{-1}$, $c\ts b' c^{-1}$ and\/ $d$ pairwise commute.
\end{enumerate}
Moreover, if\/ $\cG$ is our special group\/ ${\rm GA}^{}_2(K)$,
its centre is trivial.   
\end{lemma}
\begin{proof}
The first three assertions simply are a reformulation of
\cite[Thm.~4.5]{MKS} in our context. From \cite[Cor.~4.5]{MKS}, we
also know that $\cent(\cG) = \cent(\cA)\cap\cent(\cE)\subset\cB$,
whenever $\cA\neq\cB\neq\cE$. This is clearly the case for our special
group $\cG={\rm GA}^{}_2(K)$.  It is not difficult to verify that
$\cent(\cA)=\{1\}$, which then establishes the last claim.
\end{proof}

The next step is a complete characterisation of the Abelian subgroups
into three types, which goes back to Moldavanskii \cite{Mol}. It was
later put into a more general framework in \cite{KS}, and a complete
account is also contained in \cite[Sec.~0]{Wright}. We first rephrase
\cite[Thm.~0.3]{Wright} in our terminology, but still for a general
setting. We will then specialise step by step.
\begin{theorem} \label{abelian-subgroups}
   If\/ $\cH$ is an Abelian subgroup of\/ $\cG=\cA\amalgam\cE$, it is
   precisely of one of the following three types.
\begin{itemize}
\item[(T 1)]
      $\cH$ is conjugate to a subgroup of\/ $\cA$ or
      to a subgroup of\/ $\cE$.
\item[(T 2)] $\cH$ is not conjugate to any subgroup of\/ $\cA$ or\/
      $\cE$, but there exists a nested chain of subgroups\/ $\cH^{}_0
      \subset \cH^{}_1 \subset \dots \subset \cH^{}_i \subset \dots$
      such that\/ $\cH = \bigcup_{i=0}^{\infty} \cH_i$, where each\/
      $\cH_i$ is conjugate to a subgroup of\/ $\cB$. This chain is
      inevitably infinite and non-stationary.
\item[(T 3)] $\cH = \cF\times\langle g \rangle$, where\/ $\cF$ is
      conjugate to a subgroup of\/ $\cB$, and\/ $g$ is a CR element,
      hence not of finite order and not conjugate to any element of\/
      $\cA$ or\/ $\cE$.  \qed
\end{itemize}
\end{theorem}

As is immediate, type 2 subgroups are the more delicate ones to deal
with. We will now focus on our special group of polynomial
automorphisms $\cG=\cG^{}_K$, which admits further simplifications.
For completeness, we will consider all three types here, even though
later on we will mainly need Abelian subgroups of type 3.

\begin{prop} \label{type-one}
   Let\/ $\cH$ be an Abelian subgroup of\/ $\cG = \cA \amalgam \cE$.
   Then, the following three assertions are equivalent.
\begin{enumerate}
\item $\cH$ is of type\/ $1$.
\item $\cH$ is conjugate to a subgroup of either\/ $\cA$ or\/ $\cE$.
\item $\cH$ is of bounded length, i.e., 
      $\max\{\len(g)\mid g\in\cH\}<\infty$.
\end{enumerate}
\end{prop}
\begin{proof}
(1) $\Longleftrightarrow$ (2) is the definition. \\
(2) $\Longrightarrow$ (3):
Each element of $\cA\cup\cE$ has length $0$ or $1$, which is then
also true of any subgroup $H$ of $\cA$ or $\cE$. Since
$  \len(g h g^{-1}) \le  \len(g) + \len(h) + \len(g^{-1})$,
any conjugate subgroup $gHg^{-1}$ is then of bounded length, too. \\
(3) $\Longrightarrow$ (2):
This is proved in \cite[Prop.~0.35]{Wright}, or follows from
Fact~\ref{serre-fact}.  
\end{proof}

{}From now on, let $U = U_K$ denote the group of roots of $1$ in $K$,
so $U=\{\pm 1\}$ for $K=\RR$ and $U=\{z\in S^1 \mid z^n=1 \mbox{ for
some } n\in \NN\}$ for $K=\CC$ (with $S^1$ the unit circle in $\CC$).
Moreover, let $U(n)$ denote the (multiplicative) subgroup of $n$-th
roots of unity of $K$. Note that $U(n)$ is a finite cyclic group
\cite[Thm.IV.1.9]{Lang}, the order of which divides $n$. If $n$ is a
power of $\ch (K)$, one has $U(n)=\{1\}$, the trivial group. If $n$
is not divisible by $\ch (K)$, and if $K$ is algebraically closed, one
has $U(n)\simeq C_n$; without algebraic closure, $U(n)$ can be a
genuine subgroup of $C_n$, as happens in $\RR$ versus $\CC$, compare
\cite[Sec.~VI.3]{Lang} for more.

\smallskip
{}Following \cite[Thm.~1.21 and Cor.~1.22]{Wright}, one can summarise
the situation of type 2 subgroups of $\cG$ as follows.
\begin{prop} \label{type-two}
   The group\/ $\cG$ does not contain Abelian subgroups of type\/ $2$
   if\/ $K$ is a finite field, or if\/ $K$ has characteristic\/ $0$
   and finite\/ $U$.

   Otherwise, if\/ $\cH$ is a type\/ $2$ Abelian subgroup of\/ $\cG$,
   the necessarily non-stationary subgroup chain\/ $(\cH_i)_{i\ge 0}$
   of Theorem~$\ref{abelian-subgroups}$ satisfies one of the following
   two conditions.
\begin{enumerate}
\item Each\/ $\cH_i$ is conjugate to a finite subgroup of\/ $K^2$,
      viewed as a subgroup of\/ $\cA$.
\item Each\/ $\cH_i$ is conjugate to a subgroup of the diagonal
      matrices of the form \newline
      $\{\diag(u,u^m)\mid u\in U(n)\}$,
      where\/ $m$ and\/ $n$ are coprime integers that depend on\/ $i$.
\end{enumerate}
Moreover, if\/ $K$ has characteristic\/ $0$, only case\/ $(2)$ is 
possible.      \qed
\end{prop}

In general, all situations can occur, see Examples 2.2 and 2.5 of
\cite{Wright}. If restricting to $\QQ$, $\RR$ or $\CC$, we are in the
case of characteristic $0$. But while for $K\in\{\QQ,\RR\}$ we do not
have type 2 subgroups (since then $U_K = K\cap U_{\CC} = \{\pm 1\}$,
so that we cannot have non-stationary subgroup chains), this is not so
for $K=\CC$.

\smallskip
{}Finally, we recall \cite[Thm.~1.24]{Wright}.

\begin{prop} \label{type-three}
   If\/ $\cH$ is a type $3$ Abelian subgroup of\/ $\cG$, then\/ $\cH =
   \cF \times \langle g \rangle$, where\/ $g\in\cH$ is a CR element
   $($hence not of finite order and not conjugate to an element of\/
   $\cA$ or\/ $\cE$\hspace{1pt}$)$, and\/ $\cF$ is a subgroup of\/
   $\cH$ such that one of the following two conditions holds.
\begin{enumerate}
\item $\cF$ is conjugate to a subgroup of\/ $K^2$, the latter
      viewed as a subgroup of\/ $\cA$.
\item $\cF$ is conjugate to a subgroup of the diagonal
      matrices of the form \newline
      $\{\diag(u,u^m)\mid u\in U(n)\}$,
      where\/ $m$ and\/ $n$ are $($fixed\/$)$ coprime integers.
      \newline
      In particular, $\cF$ is a finite cyclic group.
\end{enumerate}
Once again, if\/ $K$ has characteristic\/ $0$, only case\/ $(2)$ is
possible.  \qed
\end{prop}

As follows from Examples 2.9 and 2.10 of \cite{Wright}, both
possibilities of Proposition~\ref{type-three} can be realised in
general fields. We will provide another example of this in the next
section. If $\ch(K)=0$ or if $K$ is a finite field, the group $\cF$ is
finite. Whether this is generally the case, as addressed on p.\ 613 of
\cite{Wright}, does not yet seem to have been resolved \cite{W2}.

\section{Symmetries}  \label{symmetries}

We now turn our attention to the symmetry group $\cS(f) =
\cent^{}_{\cG} (f) = \{ h\in \cG\mid fh = hf\}$ for $f \in \cG$. In
particular, we would like to know its structure, e.g., whether it is
Abelian.  Although this need not be the case in general, the knowledge
of the Abelian subgroups reviewed in the previous section will prove
most useful to determine the structure of $\cS(f)$.

We are mainly interested in the case that $f$ is a CR element, because
these are dynamically the most interesting ones.  This is also
justified by the observation that, due to Lemma~\ref{commute}, the
investigation of the symmetries of other elements can essentially be
handled within the factors $\cA$ or $\cE$. Even though this is a task
in itself (note that it also includes the analysis of point and space
groups in the plane, hence cases where $\cS(f)$ is {\em not\/}
Abelian), it is more or less decoupled from $\cG$, due to the very
structure of $\cG$ as an amalgamated free product. In fact, the
analysis of space groups is essentially restricted to $\cA$, see
\cite[Sec.\ 4.5]{CM} for details.

As to $\cE$, consider an element $e\in\cE$ of the form given in
\eqref{elementary} with $\alpha=1$ and $v=0$. It is immediate that it
always commutes with the simple translation $t\!:\, x' = x + 1, \, y'
= y$, where $\langle t \rangle$ is isomorphic with $C_{\infty}$
(resp.\ $C_p$) if $\ch(K)=0$ (resp.\ $\ch(K)=p$ with $p$ prime).  Now,
let the polynomial $P$ from $e$ be {\em odd} (i.e., $P(-y) = - P(y)$),
and consider an arbitrary field $K$ with $\ch(K)\neq 2$, so that
$-1\neq 1$.  Clearly, $e$ now also commutes with the mapping defined
by $I = \diag(-1,-1) \in {\rm GL}(2,K)$, but $t$ and $I$ do {\em
not\/} commute (one has $I\circ t = t^{-1} \circ I$).  Though both
$\langle e,t\rangle$ and $\langle e,I\rangle$ are Abelian subgroups of
$\cE$, hence Abelian subgroups of $\cG$ of type 1, $\cS(e)$ is never
Abelian in this case, as it contains $\langle t,I\rangle$, which is a
dihedral group of infinite order ($D_{\infty}$) or of order $2p$
(denoted by $D_p$).

\smallskip
We now specialise to investigate $\cS(f)$ for $f\in\cG$ a CR element.
Two themes will run through our investigation of the symmetries of
such an element: (i) we profit from studying the ``local'' symmetry
group $\langle f, h \rangle$ with $h\in\cS(f)$, which is an Abelian
subgroup of $\cS(f)$, though $\cS(f)$ itself might not be Abelian;
(ii) the order of possible symmetries (and later of reversing
symmetries) is driven by the nature of the roots of unity in the
chosen field $K$.

We start with a simple observation which highlights a first difference
with the above example of the (non-CR) element from $\cE$.
\begin{lemma} \label{exclude}
    If\/ $f\in\cG$ is a CR element of\/ $\cG$, it cannot be contained
    in any Abelian subgroup of\/ $\cG$ of type\/ $1$ or\/ $2$.
\end{lemma}
\begin{proof}
A CR element is not conjugate to any element of
$\cA$ or $\cE$, by Fact~\ref{conjugate}, hence cannot lie in an
Abelian subgroup of type 1, by Proposition~\ref{type-one}.

On the other hand, by Theorem~\ref{abelian-subgroups} and
Proposition~\ref{type-two}, all type 2 Abelian subgroups $\cH$ are
obtained as inductive limits of a sequence $(\cH_i)_{i\ge 0}$ of
nested groups, each of which is conjugate to a subgroup of
$\cB$. Since neither the CR element $f$ nor any of its (finite) powers
can be an element of any of these subgroups, $f$ is not an element of
$\cH$ either.
\end{proof}

\noindent {\sc Example}: Let $K$ be an arbitrary field, and consider 
the mapping 
\[
   f: \quad x' = y \; , \quad y' = x + Q(y)
\]
with the polynomial $Q(y) = y^p - y$, where $p$ is a prime.  This is a
CR element of length $2$, whose square would be in the generalised
standard form \eqref{defgenstan} with $P_1=P_2=Q$.  If $K$ has
characteristic $p$ (e.g., if $K = {\mathbb F}_p$), it is easy to check
that $Q(y+1) = Q(y)$, because $\binom{p}{\ell} = 0$ (mod $p$) for all
$1 < \ell < p$ (this is equivalent to the existence of the Frobenius
endomorphism in characteristic $p$, defined by $y\mapsto y^p$, cf.\
\cite[p.\ 179]{Lang}). As a consequence, $f$ commutes with the
translation $t\in \cB \!:\, x'=x+1, \, y'=y+1$, the latter generating
the cyclic group $C_p$.

If we now restrict to {\em odd\/} primes (i.e., $p\neq 2$), the
polynomial $Q$ is odd, and $f$ also commutes with $I=\diag(-1,-1)$.
The latter, in turn, does {\em not\/} commute with $t$, and $\langle
t,I\rangle \simeq D_p$.  The minimal example emerges for $p=3$, where
$D_3\simeq S_3$ (with $S_n$ the symmetric group) is the smallest
non-Abelian group.  Clearly, $\cS(f)$, which contains $\langle
t,I\rangle$, is not Abelian either.

Note, however, that both $\langle f,t\rangle$ and $\langle f,I\rangle$
are Abelian subgroups of $\cG$ of type 3, fitting cases (1) and (2) of
Proposition~\ref{type-three}, respectively.

\medskip
This example shows that, in general, the symmetry group of a CR
element will not be Abelian, but also that interesting new phenomena
occur when one works over finite fields or over fields with
characteristic $\neq 0$.

To use the knowledge of Abelian subgroups of $\cG$ of the previous
section, it seems a reasonable strategy to restrict, as far as
possible, to ``local'' symmetries, i.e., to the groups generated by a
CR element $f$ together with a single symmetry. Then,
Lemma~\ref{exclude} has the following consequence.
\begin{prop}  \label{structure}
   Let\/ $f\in\cG$ be a CR element of\/ $\cG=\cG_K$ and\/ $h$ be a
   symmetry of\/ $f$, i.e., $h\in\cS(f)$. Then, $\langle f,h \rangle$
   is an Abelian subgroup of\/ $\cG$ of type\/ $3$.

   Moreover, if\/ $\ch(K)=0$, one has\/ $\langle f,h \rangle \simeq
   C_{\ell}\times C_{\infty}$ for some\/ $\ell\in\NN$. This means that
   one either has\/ $f^k = h^m$ for some\/ $k,m\in\ZZ\setminus\{0\}$,
   or\/ $\langle h \rangle$ is a finite cyclic group.
\end{prop}
\begin{proof}
Even though $\cS(f)$ itself need not be Abelian,
each subgroup of the form $\langle f,h\rangle$ certainly
is. Since $f\in\langle f,h\rangle$ and $f$ is a CR element, the first
claim follows from Lemma~\ref{exclude}.

If $K$ is a field of characteristic $0$, part~(2) of
Proposition~\ref{type-three} says that the Abelian subgroups of type 3
are all of the form $\cF\times C_{\infty}$, where $\cF$ is isomorphic
with a subgroup of the cyclic group $C_n$, for a suitable $n$.
Consequently, $\cF\simeq C_{\ell}$ for some divisor $\ell$ of $n$.

So, we have $\langle f,h \rangle \simeq \langle t \rangle \times
\langle g \rangle$ with $t^{\ell}=1$ and $g$ an element of infinite
order.  Clearly, $f=t^{\epsilon} g^{r}$ for some $0\le\epsilon < \ell$
and $0 \ne r\in\ZZ$, hence $f^{\ell} = g^{\ell r}$. By the same
argument, $h=t^{\epsilon'}g^q$ and $h^{\ell} = g^{\ell q}$ for some
$q\in\ZZ$, possibly $0$. Consequently, $f^k = h^m$ with $k=\ell q$ and
$m=\ell r$.

If $k\neq 0$, we must also have $m\neq 0$, since $f$ is not of finite
order.  This gives the first possibility claimed, where $h$ is not of
finite order.  If $k=0$, one has $h^m=1$, whence $h$ is of finite
order. This is only possible for $q=0$, so $\langle h \rangle$ is
isomorphic with a subgroup of $C_{\ell}$ and hence cyclic.  
\end{proof}

\noindent {\sc Remark}: If $\ch(K)=0$, Proposition~\ref{structure}
excludes, for $f$ a CR element, the existence of a subgroup of
$\cS(f)$ of the form $C_{\infty} \times C_{\infty}$ that contains $f$,
i.e., the existence of an infinite order symmetry which is independent
of $f$. This observation forms the basis of a result of Veselov
\cite{V1,V2} that an (area-preserving) polynomial automorphism cannot
possess a polynomial or rational integral $I(x,y)$ that is preserved
under iteration of $f$.

\medskip
The following example illustrates (e.g., when $K=\CC$) that symmetries of
arbitrarily large finite order can indeed occur.

\medskip
\noindent {\sc Example}: Let $K$ be a field, with unit group $U_K$.
Take $f\in\cG^{}_{K}$ of the generalised standard form
\eqref{defgenstan} and look for a linear symmetry $h$ defined by the
matrix $\diag(\lambda,\mu)\in {\rm GL}(2,K)$, hence with
$\lambda,\mu\in K^*=K\setminus \{0\}$.  One finds that $f\circ h =
h\circ f$ if and only if the polynomials $P_1$ and $P_2$ of
\eqref{defgenstan} satisfy
\[
   P_1 (\mu z) \; = \; \lambda\ts P_1 (z)  \quad \mbox{and} \quad
   P_2 (\lambda z) \; = \; \mu\ts P_2 (z) \ts .
\]
In particular, unless $P_1 = P_2 = 0$ or
one polynomial vanishes while the other is a non-constant monomial,
there is no solution except when $\lambda$ and $\mu$ are roots of
unity, i.e., when $\lambda,\mu\in U$.

If $\lambda\in U(n)$ is a {\em primitive} $n$-th root of $1$ in $K$,
and if $\mu = \lambda^{-1}$, the order of $h$ is $n$.  Moreover, the
polynomial condition is satisfied if $z\ts P_1 (z)$ and $z\ts P_2 (z)$
are actually polynomials in $z^n$ without constant term.  Similarly,
if $\mu=\lambda$, one needs $\frac{1}{z} P_i(z)$ to be a polynomial in
$z^n$, for $i\in\{1,2\}$.  Consequently, if $K=\CC$, symmetries of any
finite order are possible.

\medskip
\noindent {\sc Remark}: Proposition~\ref{structure} shows that, when
$\ch(K)=0$, a symmetry $h$ of a CR element $f$ is:\ either (i) of
finite order, conjugate to a diagonal matrix with entries from the
roots of unity; or (ii) of infinite order and $h$ and $f$ are both
roots of a common CR element. Furthermore, from the proof of
Proposition~\ref{structure} and part~(2) of
Proposition~\ref{type-three}, it follows immediately that $f=(c\,
\diag(u,u^m)\, c^{-1})^{\epsilon} (d\, G\, d^{-1})^{r}$ for some $u\in
U(\ell)$, $0\le\epsilon < \ell$, $c,d,G \in \cG$ with $G$ having a
CRNF.  Moreover, the two bracketed elements commute. Equivalently,
$f$ is conjugate to $(\diag(u,u^m))^{\epsilon}\, g^{r}$ where $g$ is
CR and commutes with the linear map defined by the diagonal
matrix. {}For $K=\CC$, Theorem 1 (or, in more detail, Theorem 7 and
Corollary 9) of \cite{GM2} is a stronger result in this spirit,
obtained by constructive means.  It is shown that $f$ is conjugate to
$(\diag(u,u^m))^{\epsilon}\, H^{r}$, with $H=h_m \circ \ldots \circ
h_1$, where the H\'enon maps $h_i$ are defined by
\begin{equation} \label{defgenhenonmap} 
 h_i \; : \quad   x' \; = \; y  \; , \quad y' \; = \; -\delta_i 
     \ts x + Q_i(y)\ts .  
\end{equation}
It follows from \cite{FM} that every CR element of $\cG$ is conjugate
to a composition $h_m \circ \ldots \circ h_1$ for some $m \ge 1$
(\cite{FM} also shows that some normalisation can be made to each
$Q_i(y)$).  With the choice of coset representatives $\cI$ of
$\eqref{I-def}$ and $\cJ$ of $\eqref{J-def}$, leading to the normal
form \eqref{product-form}, we have a similar result:\ namely, every CR
element is conjugate to a uniquely-expressed composition $(a_m \circ
e_m) \circ \ldots \circ (a_1 \circ e_1)$, resp.\ one with an extra
$b\in\cB$ in front of it.  Note that
\begin{equation} \label{defaiei} 
 a_i \circ e_i \; : \quad   x' \; = \; y  \; , \quad y' \; = \; 
     \ts x + (y^2 P^{}_i (y) + \beta^{}_i y)\ts   
\end{equation}
is an orientation-reversing H\'enon map.

\medskip
In the cases $K=\QQ$ and $K=\RR$, there are more severe restrictions on
the nature of symmetries and, in fact, $\cS(f)$ turns out to be Abelian.
\begin{theorem}  \label{involutions}
   Let\/ $K$ be a field of characteristic\/ $0$ with group of roots of
   unity\/ $U\simeq C^{}_2$ $($which includes the cases\/ $K=\QQ$
   and\/ $K=\RR$\hspace{1pt}$)$, and let\/ $f$ be a CR element of\/
   $\cG$. Then, any symmetry of\/ $f$ in\/ $\cG$ of finite order must
   be the identity or an involution.

   Moreover, the symmetry group of\/ $f$ in\/ $\cG$ can contain at
   most one nontrivial finite group, which is then of the form\/
   $\langle s \rangle$ with\/ $s$ an involution that is conjugate to\/
   $I = \diag(-1,-1)$.
\end{theorem}
\begin{proof}
Let $h$ be any element of $\cS(f)$. Due to the
assumptions, $\langle f,h \rangle$ is an Abelian subgroup of type 3,
hence equals $\cF\times\langle g \rangle$ for some CR element $g$ and
some {\em finite\/} group $\cF$. Since we are in case (2) of
Proposition~\ref{type-three}, $\cF$ is isomorphic to a subgroup of the
group $U$ of roots of unity in $K$, hence to the trivial group or
$C^{}_2$. So, $\cF=\langle s \rangle$ with $s^2 = 1$. Since $\cF$
contains all elements of $\langle f,h \rangle$ of finite order, and
$h$ was an arbitrary symmetry, this shows that any symmetry of $f$ of
finite order must be $1$ or an involution.

Clearly, $f$ itself is in $\cS(f)$, but it is a CR element, hence not
of finite order. So, we must have $f = s^{\epsilon}g^m$ for $\epsilon$
either $0$ or $1$ and some nonzero integer $m$, hence $f^2 =
g^{2m}$. So far, we have established that $f$, together with any
single symmetry of it, generates an Abelian group of the form
$\cF\times C_{\infty}$ with $\cF$ the trivial group or $C^{}_2$. We
now need to understand better how different groups of this kind fit
together as subgroups of $\cS(f)$.
 
So, let us assume that $\cS(f)$ contains two different involutions,
$s^{}_1$ and $s^{}_2$ say. Then, also the product $s = s^{}_1 s^{}_2$
commutes with $f$. Our previous argument applies to $h=s$, so we have
$\langle f,s \rangle = \cF\times\langle g \rangle$ with $\cF$ the
trivial group or $C^{}_2$, and $g$ a CR element with $f^2=g^{2m}$ for
some nonzero $m\in\ZZ$. So, $\cF=\langle t \rangle$ with $t^2=1$,
hence $s=t^{\epsilon} g^k$ with $\epsilon\in\{0,1\}$ and $k\in\ZZ$.

If $s$ is not of finite order, one has $k\neq 0$ and $s^{2m} =
(t^{\epsilon} g^k)^{2m} = g^{2mk} = f^{2k}$. With $s = s^{}_1 s^{}_2$,
observe $s_i s^{2m} s_i = s^{-2m}$, for $i\in\{1,2\}$, so that
$s^{}_1$ and $s^{}_2$ are reversors of $f^{2k}$. But the $s_i$ are
also symmetries of $f$, hence of $f^{2k}$, and we obtain $f^{-2k} =
s^{}_i f^{2k} s^{-1}_i = f^{2k}$ which would imply $f^{4k}=1$ --- a
contradiction.

So, $s=t^{\epsilon} g^k$ must be of finite order. This implies $k=0$
(because $g^0=1$ is the only finite order element of $\langle g
\rangle \simeq C_{\infty}$), hence $s=t^{\epsilon}$ and $s^2=1$. Since
$s^{}_1 \neq s^{}_2$ by assumption, we know that $s\neq 1$, and $s$
must be an involution.  This also implies that $s^{}_1$ and $s^{}_2$
commute.  But $s^{}_1 \neq s^{}_2$ now means that we have an Abelian
subgroup $\langle s^{}_1 \rangle \times \langle s^{}_2 \rangle \times
\langle f \rangle \simeq C_2 \times C_2 \times C_{\infty}$ of $\cG$
which must be of type 3. However, the finite group $\cF$ here is
Klein's 4-group, which is not cyclic. This contradicts
Proposition~\ref{type-three}.

Consequently, there can be at most {\em one\/} true involution which
commutes with $f$, which shows the claim about the finite subgroup of
$\cS(f)$. In fact, part~(2) of Proposition~\ref{type-three} implies
that $s$ is conjugate to $I=\diag(-1,-1)$.  
\end{proof}

Let us draw some further conclusions from Theorem~\ref{involutions},
under the assumptions given there.  If $g$ is an element with CRNF and
$h\in\cG$ satisfies $g^n=h^n$ for some positive integer $n$, also $h$
must have CRNF (this follows from a simple argument involving the
length of the elements and their powers and the fact that the power of
a CR element $h$, after reduction to normal form, must start and end
with elements of the same type, i.e., from $\cI$ or $\cJ$, as $h$
itself). In fact, the only possibility is $h=bg$ for some
$b\in\cB$. Clearly, $h$ commutes with $g^n$. This implies that
$b\in\cS(g^n)$, and $b$ must be of finite order (since otherwise
$\langle b,g^n\rangle \simeq C_\infty \times C_\infty$, which is
impossible). By Theorem~\ref{involutions}, either $b=1$ (whence
$b\in\cS(g)$) or $b$ is the unique involution in $\cS(g^n)$. In the
latter case, also $gbg^{-1}$ is an involution in $\cS(g^n)$, hence
$gbg^{-1} = b$ by uniqueness, and $b\in \cS(g)$. Since this applies to
general CR elements by conjugacy, we have
\begin{fact} \label{roots}
  If\/ $K$ is a field with\/ $\ch(K)=0$ and\/ $U_K\simeq C_2$, a CR
  element\/ $f\in\cG$ has at most one $n$-th root in\/ $\cG$ for\/
  $n$ odd, and at most two for\/ $n>0$ even. If two roots exist, one
  is obtained from the other by multiplication with the unique
  involution in\/ $\cS(f)$.  \qed
\end{fact}

\begin{coro} \label{realisAbelian}
   Let the assumptions be as in Theorem~$\ref{involutions}$, with\/
   $f\in\cG$ a CR element. Then, $\cS(f)\simeq \cF\times C_{\infty}$,
   where\/ $\cF$ is either the trivial group or\/ $C_2$, and\/
   $C_{\infty}$ is generated by a CR element. In particular, $\cS(f)$
   is Abelian.
\end{coro}
\begin{proof}
If $\cS(f)$ contains any nontrivial element of finite order at all,
$s$ say, it must be an involution and is unique, due to
Theorem~\ref{involutions}.  For an arbitrary $h\in\cS(f)$, also
$hsh^{-1}$ is an involution, hence $hsh^{-1}=s$. So, $s$ commutes with
all elements of $\cS(f)$ and is thus an element of its centre.
Moreover, $s$ is conjugate to $I$ by Theorem~\ref{involutions}.

No element of $\cS(f)$ other than $1$ and possibly $s$ can be of
finite order. In fact, they must all be CR elements (otherwise, we
would obtain an Abelian subgroup of the form $C_{\infty}\times
C_{\infty}$, which is impossible). If $g\neq 1$ is such an element, we
know from part (3) of Lemma~\ref{commute} that $f = cbc^{-1}\, d^k$
and $g = cb'c^{-1}\, d^{\ell}$ with $b,b'\in\cB$, $d\in\cG$, and
suitable $k,\ell\in\ZZ$. Also, $cbc^{-1}$, $cb'c^{-1}$ and $d$
pairwise commute, so must all be elements of $\cS(f)$.  Consequently,
each of $b$ and $b'$ can only be $1$ or conjugate to $s$, while
$d$ must be a CR element (and $k,\ell\neq 0$).

Let us now, without loss of generality, assume that $f$ has CRNF, so
$\len(f) = 2n$ with $n\ge 1$.  This implies that the equation $f=h^m$,
with $h\in\cG$, can at most have a solution if $m$ divides $n$ and if
$h$ is another element with CRNF. Clearly, $h$ itself commutes with
$f$.  An analogous restriction applies to the equation $f= s h^m$ with
$s$ an involution from $\cS(f)$, because then $sf = h^m$, and $h$
commutes with $sf=fs$.

In both cases, we can invoke Lemma~\ref{commute} once more. Since $f$,
by Fact~\ref{roots}, can only have one odd root and at most two even
roots, there must be a fundamental element $h$ which, possibly
together with the unique involution $s$, can be used for {\em all\/}
symmetries $g$ of infinite order, so that $g=s^{\epsilon} h^m$ for
some $\epsilon\in\{0,1\}$ and some $m\in\ZZ \setminus \{0\}$, with $m$
even if $\epsilon=1$. This shows that $\cS(f) = \cF\times \langle
h\rangle$, where $\cF$ is the trivial group or $C_2$, and $\langle h
\rangle\simeq C_{\infty}$.
\end{proof}

\noindent {\sc Remark}:
It would be interesting to know whether $\cS(f)$ is always Abelian
for the case $K=\CC$.

\smallskip
In view of Theorem~\ref{involutions}, and also for later use as
potential reversors, it is of particular interest to know the
involutions in $\cG$, up to conjugacy. Since for $\ch (K) = 2$ one has
$1=-1$, so that there are no $2^k$-th roots of unity except $1$,
usually no involutions (or elements of order $2^k$) exist in ${\rm
GL}(2,K)$ that are of interest to us here (though new involutions in
$\cG$ will show up, such as the elements of $\cJ$).  Consequently, we
will exclude fields of characteristic $2$ in what follows.
\begin{lemma} \label{fin-ord}
    If\/ $K$ is a field with\/ $\ch(K)\neq 2$, the possible 
    involutions in\/ $\cE$ are
\begin{equation} \label{case-1}
   e\! : \; x' \, = \, -x + P(y) \, , \; 
   y' \, = \, y
\end{equation}
   with arbitrary polynomial\/ $P\in K[y]$, or
\begin{equation} \label{case-2}
   e\! : \; x' \, = \, \alpha\, x + P(y) \, , \; 
   y' \, = \, -y + v
\end{equation}
  with arbitrary\/ $v\in K$, $\alpha\in\{\pm 1\}$ and a polynomial\/
  $P\in K[y]$ that satisfies\/ $P(v-y) = - \alpha\, P(y)$.

   Moreover, if\/ $K$ is a field of characteristic\/ $0$ with\/ $U
   \simeq C^{}_2$, any element of\/ $\cE$ of finite order is either
   the identity or an involution.
\end{lemma}
\begin{proof}
Consider $e\in\cE$, parametrised as in (\ref{elementary}). Then, the
first claim is a straightforward calculation around the equation
$e^2=1$.

{}For the second claim, write $e^n$, for integer $n\ge 0$, as $x\mapsto
x_n$ and $y\mapsto y_n$.  Setting $v^{}_0=0$ and $v_n =
(1+\beta+\ldots+\beta^{n-1})v$ for $n\ge 1$, one has $y_n = \beta^n y
+ v_n$, and a direct calculation gives
\[
   x_n \; = \; \alpha^n x + \sum_{\ell=0}^{n-1}
   \alpha^{n-1-\ell} P(y^{}_\ell) \, .
\]
Clearly, $e^n=1$ implies $\alpha^n=\beta^n=1$, $v_n=0$ and
$\sum_{\ell=0}^{n-1} \alpha^{n-1-\ell} P(y^{}_\ell) = 0$.  In
particular, $n$ odd is impossible due to $U\simeq C_2$, unless 
$\alpha =\beta =1$.

The case $\beta=1$ means $v_n=n v$, hence $v=0$ because $\ch(K)=0$. If
also $\alpha=1$, the polynomial must be $P=0$, and hence $e=1$. If
$\alpha=-1$, $e^n=1$ is true for all even $n$ and arbitrary $P$, but
one actually has $e^2=1$.

{}For the case $\beta=-1$, one has $v_n=0$ for all even $n$, and $e^n=1$
follows for all $P\in K[y]$ with $\alpha\, P(y) + P(v-y) = 0$. 
Clearly, one has $e^2=1$ in these cases, too.  
\end{proof}

This classifies the involutions in $\cE$ for $\ch(K)\neq 2$.  The
involutions in $\cA$ are the elements of the form $(\bs{a},M)\neq
(\bs{0},1)$ that satisfy
\begin{equation} \label{involA}
M^2=1 \quad \mbox{and} \quad M\bs{a}=-\bs{a}. 
\end{equation}
Investigating the first of these requirements, we find
\begin{lemma} \label{affine-involutions}
   Let\/ $K$ be an arbitrary field with\/ $\ch(K)\neq 2$. If\/ $M$ is
   an involution in\/ ${\rm GL}(2,K)$, it is either\/
   $I=\diag(-1,-1)$, or it is ${\rm GL}(2,K)$-conjugate to\/
   $T=\left(\begin{smallmatrix} 0 & 1 \\ 1 & 0
   \end{smallmatrix}\right)$ or, equivalently, to\/
   $S=\diag(-1,1)$. Moreover, if an involution\/ $M\neq I$
   is not upper triangular, it is conjugate to\/ $T$
   by a matrix\/ $A\in\cT$, the subgroup of\/ ${\rm GL}(2,K)$
   of invertible upper triangular matrices.
\end{lemma}
\begin{proof}
Consider the equation $M^2=1$ with $M\in{\rm GL}(2,K)$, which is
clearly solved by $I$. An easy direct calculation shows that all other
solutions satisfy $\trace(M)=0$ and $\det(M)=-1$, hence share the
characteristic polynomial $P(x)=x^2-1=(x+1)(x-1)$. As this is then
also the minimal polynomial (we excluded $\ch(K)=2$), they also share
all polynomial invariants and must have the same rational canonical
form, compare \cite[Ch.~4.4]{AW}. Consequently, they are all ${\rm
GL}(2,K)$-similar to the Frobenius companion matrix of $P(x)$, which
is $T$. In particular, $S$ is conjugate to $T$ in ${\rm GL}(2,K)$.

If $M\neq I$ is an involution, we know that
$M=\left(\begin{smallmatrix} a & b \\ c & -a \end{smallmatrix}\right)$
with $a^2+bc=1$. If $M$ is {\em not\/} upper triangular, we also know
that $c\neq 0$. Then, it is easy to check that
$A=\left(\begin{smallmatrix} c & -a \\ 0 & 1
\end{smallmatrix}\right)$, which lies in $\cT$, satisfies
$AMA^{-1}=T$.  
\end{proof}

Continuing the investigation of affine involutions, but
now also considering how the involutions of
Lemma~\ref{fin-ord} are related to linear ones, we obtain
(compare also \cite[p.\ 57]{Kraft})
\begin{prop} \label{elementary-to-linear}
   If\/ $K$ is a field with\/ $\ch(K)\neq 2$, all involutions in\/
   $\cA\cup\cE$, and hence in $\cG$, are conjugate to linear maps.
   More concretely, they are conjugate to either\/ $I=\diag(-1,-1)$ or
   to\/ $S=\diag(-1,1)$, equivalently to\/
   $T=\left(\begin{smallmatrix} 0 & 1 \\ 1 & 0
   \end{smallmatrix}\right)$, where\/ $T\in\cI$ of\/ $\eqref{I-def}$.
\end{prop}
\begin{proof}
First, consider the affine involutions which must satisfy
\eqref{involA}. From Lemma~\ref{affine-involutions}, one possibility
based around $M=I$ is $g\in\cB$ defined by $x'=-x+u,\,y'=-y+v,\,$ with
arbitrary $u$, $v\in K$. Taking $h\in\cB$ via $x' = x-u/2$ and $y'=
y-v/2$, one finds $h g h^{-1}$ is the linear map defined by $I$. The
second possibility for affine involutions, from
Lemma~\ref{affine-involutions}, consists of those that are ${\rm
GL}(2,K)$-conjugate to $g\in\cA$ defined by $x'=y+u$ and $y'=x-u$
(noting from \eqref{involA} that the entries of $\bs{a}$ must have
opposite signs when $M=T$). However, this $g$ is itself conjugate in
$\cB$, via the above-mentioned $h$ with $v = -u$, to the linear map
defined by $T$. Finally, it is clear that $T$ is conjugate in ${\rm
GL}(2,K)$ to the matrix $S$.

We turn now to the elementary involutions as described in
Lemma~\ref{fin-ord}. Consider the involution $e\in\cE$ from
\eqref{case-1}.  Defining $h\in\cE$ by $x' = x-P(y)/2$ and $y'=y$, one
can easily check that $h e h^{-1}$ is the linear map defined by
$S$. This establishes a conjugacy within $\cE$.

Next, consider $e\in\cE$ from \eqref{case-2}, with $\alpha=\{\pm 1\}$,
$v\in K$ and $P(v-y)=-\alpha\ts P(y)$. Define $h\in\cG$ via $x' =
y-v/2$ and $y'=x+P(y)/2\alpha$ (which has the inverse $h^{-1}$ given
by $y'=x+v/2$ and $x'=y-P(y')/2\alpha$). A short calculation using the
symmetry property of $P$ confirms that $h e h^{-1}$ is the linear map
defined by the matrix $\diag(-1,\alpha)$, hence either $I$ or $S$.
\end{proof}

Comparing Proposition~\ref{elementary-to-linear} and
Lemma~\ref{affine-involutions}, it is worth pointing out that it is
sometimes useful in deriving normal forms to have the freedom to use
either of the ${\rm GL}(2,K)$-conjugate matrices $S$ or $T$, where $T$
is an element of $\cI$, while $S$ is not. This is particularly true
when we study reversing symmetries in the next section
(cf. Proposition~\ref{invnormform}).  In the case of symmetries, $I$
is the important involution, as will turn out shortly.

It will prove useful to define the so-called {\em poly-degree}\/ of an
element $g\in\cG\setminus\cA$.  If $g$ is in normal form
(\ref{product-form}), the poly-degree is defined by
\begin{equation} \label{pol-deg-def}
   \pd(g) \; = \; 
   \big(\deg(e^{}_m),\dots,\deg(e^{}_2),\deg(e^{}_1)\big),
\end{equation}
where we drop the last entry if $e^{}_1$ is missing in the normal form. 

\begin{theorem} \label{symnormform}
   Let $K$ be a field with $\ch(K)\neq 2$ and let\/ $f\in\cG$ be a CR
   element. If $f$ has a symmetry in $\cG$ which is an involution,
   this symmetry is conjugate to $I=\diag(-1,-1)$. Also, $f$ is
   conjugate to an element with CRNF
\begin{equation} \label{involform}
   f' \; = \; b\circ a^{}_m \circ e^{}_m \circ\ldots\circ a^{}_1 
              \circ e^{}_1 \ts .
\end{equation}
    In the expression\/ $\eqref{involform}$, $m \ge 1$, $a^{}_i \in
    \cI$ of\/ $\eqref{I-def}$, $e_i \in \cJ$ of\/ $\eqref{J-def}$ must
    have\/ $P_i(y)$ odd, $e^{}_1$ and $a^{}_m$ must appear, and\/ $b$
    is linear of the form\/ $b = \left(\begin{smallmatrix} \alpha &
    \gamma \\ 0 & \beta
    \end{smallmatrix}\right)$.
    It follows that, when $f$ has such an involutory symmetry, it is
    conjugate to a cyclically reduced element $f'$ which fixes the
    origin and has \/ $\pd(f')= (n^{}_m,\ldots,n^{}_1)$, where all\/
    $n_{i}$ are odd integers $\ge 3$.

   If $K$ is a field of characteristic\/ $0$ with group of roots of
   unity\/ $U\simeq C^{}_2$ $($which includes the cases\/ $K=\QQ$
   and\/ $K=\RR$\hspace{1pt}$)$, this gives, up to conjugacy, the
   description of all finite order symmetries and the corresponding
   normal form of\/ $f$.
\end{theorem}
\begin{proof}
Suppose $f$ has an involutory symmetry.  With
Proposition~\ref{elementary-to-linear}, we can write $f \ts (h^{-1} i
h) = (h^{-1} i h)\ts f$, where $h\in\cG$ and $i$ is the linear map
defined by $I$ or $S$.  Consequently, we have $(h f h^{-1})\, i = i\,
(h f h^{-1})$, so that a conjugate of $f$, necessarily also CR,
commutes with $I$ or $S$.
Since $\ch(K)\neq 2$, the equation $2k=0$ has only the trivial solution 
in $K$, so that $K^2$ cannot contain an involution.  Thus, we are in
the situation of case (2) of Proposition~\ref{type-three}.  Consider
$M=\diag(u,u^m)$ with $u\in U(n)$ and $n,m$ coprime.  $M$ can only be
an involution if $u=-1$ which implies that $n$ must be even. Then, $m$
must be odd, and $M=I$ is the only possibility, while $S$ is ruled out
-- a result that can also be obtained by some lengthy explicit
calculations with the normal forms.

So, let us characterise those CR elements $f'$ that commute with $I$,
equivalently those that satisfy $I\ts f'\ts I = f'$. We take for $f'$
an expression of the form (\ref{product-form}) and observe that $I$
commutes with elements $a_i\in\cI$ of (\ref{I-def}), whereas for $e_i
\in \cJ$ of (\ref{J-def}) we have $e_i \ts I = I \ts e'_i$, with
$e'_i$ obtained from $e_i$ by the replacement $P_i(y) \rightarrow
-P_i(-y)$. Note that $e'_i$ is still an element of $\cJ$. Also, $b':=
I\ts b\ts I$ is still an element of $\cB$.

The uniqueness of the normal form (\ref{product-form}) for $f'$
applied to $I\ts f'\ts I = f'$ forces $e'_i = e_i$ and $b'=b$, hence
the odd degree constraint $P_i(y)=-P_i(-y)$ in $e_i$ together with
$I\ts b\ts I = b$.  The latter implies that $b$ is linear, so $f'$
fixes the origin.  If the normal form for $f'$ so found is cyclically
reduced, at least one $e_i$ and one $a_i$ must be present by
definition. Certainly, it can be brought to the form
\eqref{involform}, possibly after a further conjugation by an element
of $\cI$.  This conjugation leaves the symmetry $I$ unchanged, so the
leading basic element of the new normal form can remain linear.  If
$f'$ is not already cyclically reduced, further conjugations by
$a_i$'s and by $e_i$'s with odd $P_i(y)$ can be used to obtain
\eqref{involform}. Again, these will leave the symmetry as $I$ because
they both commute with it. Thus, these additional conjugations, if
required, will preserve the linear nature of the leading basic element
and the oddness of the polynomials in the elementary coset
representatives.

The last statement of the theorem is simply a reminder from
Theorem~\ref{involutions} of the stronger statement that can be made
under these circumstances.  
\end{proof}

\noindent {\sc Remark}: If $f$ is a CR element, but not cyclically
reduced to begin with, a cyclically reduced element $\tilde{f}$
conjugate to $f$ can always be found in an algorithmic fashion. The
poly-degree of any such element can be used to check the necessary
condition given above on the odd entries in $\pd(f')$.  This follows
since $\pd(f')$ must be the same, up to a cyclic permutation, as the
poly-degree of the cyclically reduced element $\tilde{f}$ (from part
(3) of Proposition~\ref{conjugacy}).  In \cite{RBstandard}, as an
illustration of Theorem~\ref{symnormform}, we showed by explicit
calculation that the CR elements $f\in\cG^{}_{\RR}$ of the generalised
standard form \eqref{defgenstan} could only have symmetries of finite
order conjugate to $I$. This occurred when both $P_1$ and $P_2$ were
odd.

\medskip
However, even if a cyclically reduced element that $f$ is conjugate to
satisfies the above poly-degree requirement, a further decisive test
for an involutory symmetry still follows from
Theorem~\ref{symnormform} together with Proposition~\ref{conjugacy}.

\begin{coro} \label{involtest}
   Let $K$ be a field with $\ch(K)\neq 2$ and let $f\in\cG$ be a CR
   element. Then, $f$ has a symmetry that is an involution iff any
   cyclically reduced word to which $f$ is conjugate commutes with
   $x'=-x+u$, $y'=-y+v$, with some $u$, $v\in K$. If this cyclically
   reduced word corresponds to \eqref{defmap}, this commutation means
   $P$ and $Q$ satisfy $P(-x+u,-y+v)+P(x,y)=u$ and
   $Q(-x+u,-y+v)+Q(x,y)=v$.
\end{coro}
\begin{proof}
Let $\tilde{f}$ be a cyclically reduced word with $f=h \tilde{f}
h^{-1}$ and let $f$ have an involutory symmetry (take $h=1$ if $f$ is
already cyclically reduced). From Theorem~\ref{symnormform}, we
also know that $f$ is conjugate to a cyclically reduced word in normal
form, i.e., $f'$ of \eqref{involform}, and that $f'$ commutes with
$I=\diag(-1,-1)$ by construction.  It follows that $\tilde{f}$ and
$f'$ are two cyclically reduced words that are themselves
conjugate. By Proposition~\ref{conjugacy}, $\tilde{f}$ differs by a
cyclic permutation of the elements of $f'$, followed by conjugation
with a basic element \eqref{basic}. The cyclic permutation is itself a
conjugacy by elements $a_i$ and $e_i$ of \eqref{involform}.  It
follows that $\tilde{f}$ commutes with a conjugate of $I$, indeed the
same conjugacy used to derive $\tilde{f}$ from $f'$. As $a_i$ and
$e_i$ of \eqref{involform} commute with $I$, the only conjugacy
that can alter the symmetry of $\tilde{f}$ away from $I$ is the one by
a basic element. One easily checks, for $b$ in the form \eqref{basic},
that $b I b^{-1}$ differs from $I$ by at most a translation.  The last
statement of the result follows from forcing the form \eqref{defmap}
to commute with such an involution.  
\end{proof}

\noindent {\sc Remark}: The previous result shows that, when one deals
with a cyclically reduced element $\tilde{f}$ of $\cG$, the presence
or absence of an involutory symmetry is, in some sense, obvious.  If
present, it must be of a very simple linear (or affine) form.
Inspecting the phase portrait for the case $K=\RR$ of $\tilde{f}$, one
must see the invariance by a rotation through $\pi$
around some fixed point as a prerequisite for the existence of any
finite order symmetry (the necessity of the existence of a common
unique fixed point of both $\tilde{f}$ and the possible involutory
symmetry, if present, follows from Theorem~\ref{symnormform}).

\medskip
Another useful result, which we will need later, concerns the conjugacy 
of linear maps within the group $\cG$.
\begin{lemma}  \label{conjugacy-is-linear}
   Let\/ $f,g\in\cG$ be linear maps, defined by the matrices\/
   $A_f,A_g\in {\rm GL}(2,K)$. If\/ $f=hgh^{-1}$ for some\/ $h\in\cG$,
   then $A_f$ and $A_g$ are already conjugate within\/ ${\rm
   GL}(2,K)$.
\end{lemma}
\begin{proof}
Observe first that $\big(\dd h(\bs{a})\big)^{-1} = \dd
h^{-1}(h(\bs{a}))$, for arbitrary $\bs{a}\in K^2$, which follows from
the chain rule applied to $h^{-1} h = 1$. Since $\dd f \equiv A_f$ and
$\dd g \equiv A_g$, one then derives from differentiating $f=hgh^{-1}$
at the point $\bs{a} = h(\bs{0})$ that
\[
   A_f \; = \; \dd h(\bs{0}) A_g \big(\dd h(\bs{0})\big)^{-1}
\]
where $\dd h(\bs{0})$ clearly is an element of ${\rm GL}(2,K)$.
\end{proof}

\noindent {\sc Remark}: 
We made use of the formal differentiation
rules for polynomials here.  If one is in a setting where
diffeomorphisms are well defined, the claim can be extended
accordingly.

\section{Reversing symmetries}  \label{revsymmsec}

Recall that we denote the reversing symmetry group of an element
$f\in\cG$ by
\begin{equation*}
   \cR(f) \; = \; \{ h\in\cG\mid h f h^{-1} = f^{\pm 1} \}.
\end{equation*}
This group contains the symmetry group $\cS(f)$ as a normal subgroup, and 
the factor group $\cR(f) / \cS(f)$ is either the trivial group or $C_2$. 
In general, it is difficult to
determine these groups explicitly, but if one is in a group theoretic
setting (as we are), one can at least determine the {\em structure} of
the reversing symmetry group to some extent. This, of course, need
only be done up to conjugacy, because $\cR(hfh^{-1}) = h \cR(f)
h^{-1}$.  As before, we shall focus on elements $f\in\cG$ of infinite
order, and on CR elements in particular. This means that it actually
suffices to look at elements that possess a cyclically reduced normal
form (CRNF).

Let us start with a general observation, which is a rather direct
consequence of a result of Goodson, see \cite[Prop.~2]{Goodson} and
the generalisation mentioned afterwards. We use the general group
theoretic setting mentioned in the Introduction.
\begin{lemma} \label{rev-order}
  Let\/ $f$ be an element of infinite order, and assume that\/
  $\cS(f)= \cF \times \langle g \rangle$ where\/ $\cF$ is some finite
  group of order\/ $N$ $($not necessarily Abelian$\ts)$, and\/ $g$ is
  some generator $($then necessarily of infinite order$\ts)$.  If\/
  $r$ is a reversor of\/ $f$, then\/ $r$ is an element of\/ {\em
  finite} order. Its order is even and divides\/ $2N$.
\end{lemma}
\begin{proof}
If $r$ is a reversor, $r^2$ is a symmetry, hence $r^2=s g^m$, for some
$s\in\cF$ and some integer $m$. Note that, due to the assumption of
the direct product structure, we always have $s g = g s$, even if
$\cF$ itself is not Abelian. Since the group $\cF$ is finite and of
order $N$, we know that $s^n=1$ for some $n\neq 0$ that divides
$N$. Clearly, we then have $r^{2n}=g^{mn}$.

As $f$ is not of finite order, but clearly an element of $\cS(f)$, we
may assume $f^N=g^k$ for some (positive) integer $k$ without loss of
generality, modifying the argument just used (in particular, $k\neq
0$, while $k>0$ might require to replace $g$ by $g^{-1}$).

Since $rf = f^{-1}r$ by assumption (hence also $rf^{\ell} =
f^{-\ell}r$, for all $\ell\in\ZZ$), we choose $\ell = mnN$ and obtain
$r g^{kmn} = g^{-kmn}r$. Since $g^{kmn}=r^{2nk}$, this implies $r\,
r^{2nk} = r^{-2nk} r$ and thus $r^{4nk}=1$, i.e., $r$ is of finite
order. Since $r^{2n}=g^{mn}$, this is only possible for $mn=0$, hence
$m=0$. This implies $r^{2n}=1$, so the order of $r$ divides $2N$. If
$f$ is not of finite order, it is not an involution, and $r$ can then
not be of odd order \cite[Prop.\ 5]{lamb} (hence also $r\neq 1$).
\end{proof}

\begin{theorem}  \label{rev-order-real}
   Let\/ $K$ be a field of characteristic\/ $0$ or a finite field,
   and let\/ $f$ be a reversible CR element of\/ $\cG$,
   with reversor\/ $r$. Then, $r$ is an element of finite even order.

   If\/ $\ch(K)=0$ and if, in addition, the roots of unity in\/ $K$ are\/
   $U = \{\pm 1\}\simeq C^{}_2$, the reversor\/ $r$ is
   an involution or an element of order\/ $4$.
\end{theorem}
\begin{proof}
If $r$ is a reversor, $r^2$ is a symmetry, so $r^2\in\cS(f)$.
Consider the group $\langle f,r^2\rangle$ which is Abelian, hence of
type 3 in this case. Consequently, $\langle f,r^2\rangle =
\cF\times\langle g \rangle$ with $\langle g \rangle \simeq C_{\infty}$
and $\cF$ a finite group, by Proposition~\ref{type-three}.
Considering this as a ``local'' symmetry group of $f$, within $\langle
f,r \rangle$ say, we can invoke Lemma~\ref{rev-order} and conclude
that $r$ must be of finite even order.

If the additional assumptions on $K$ are satisfied, the finite group
$\cF$ is the trivial group or $C_2$, and we can use
Lemma~\ref{rev-order} to see that $r^4=1$. Since $r\neq 1$, it must be
an involution or an element of order $4$.  
\end{proof}

{}For fields $K$ with suitable unit group $U_K$, reversors of arbitrary
even order $\ge 2$ may exist, as the following calculation
illustrates.

\medskip
\noindent {\sc Example}: Consider a field $K$ with $\ch (K)\neq 2$ and
unit group $U_K$.  Take once again $f\in\cG^{}_{K}$ of the generalised
standard form \eqref{defgenstan} and look for a linear reversing
symmetry $r$ as defined by the matrix $\left(\begin{smallmatrix} 0 &
\mu \\ \nu & 0 \end{smallmatrix} \right)\in {\rm GL}(2,K)$. Its square
is $\diag(\lambda,\lambda)$ with $\lambda = \mu\nu$, where we assume
that $\lambda\in U_K$.  One finds that $r\circ f = f^{-1}\circ r$ if
and only if
\[
    P_1 (\nu\ts z) \; = \; -\mu\ts P_2 (z)   \quad \mbox{and} \quad
    P_2 (\mu\ts z) \; = \; -\nu\ts P_1 (z)\ts ,
\]
which also implies that $P_i (\lambda\ts z) = \lambda\ts P_i(z)$
for $i\in\{1,2\}$. Consequently, $r^2$
is a symmetry of the kind explained in the example preceding
Theorem~\ref{involutions}.

The nontrivial solutions once again occur for $\lambda$ a primitive
$n$-th root of unity, for some $n\in\NN$, provided such a $\lambda$
exists in $U_K$. For $K=\CC$, solutions exist for all $n\in\NN$.  In
these cases, the order of $r$ is $2 n$. If $\nu = -\mu$, the
polynomial condition is satisfied if $P_1 (z) = P_2 (z) = z\ts
Q(z^n)$, with some polynomial $Q$.

\medskip
\noindent {\sc Remark}: For the case $K=\CC$, reference \cite{GM2}
contains a comprehensive treatment of reversors of even order, with
illustrative examples. Also, \cite[Theorem 11]{GM2} gives
a constructive proof of Theorem~\ref{rev-order-real} above for $K=\RR$.

\medskip
Theorem~\ref{rev-order-real} motivates the benefit of knowing what
possibilities there are for elements of order $2$ and $4$ in our group
$\cG$. We have discussed the situation of elements of order $2$ in
Section~\ref{symmetries}, which we will use once more below.  Let us
now look into the remaining case when $f$ is reversible with a
reversor $r$ of order $4$ (note that we do not necessarily require
$\ch(K)=0$, although it provides an obvious motivation for this case).

\begin{theorem} \label{natureorder4}
   Let $K$ be a field with $\ch(K)\neq 2$, with a unit group\/ $U$
   that contains\/ $\{\pm 1\}$, but no primitive $4$-th root of unity
   $($thus including the case\/ $U\simeq C_2$\/$)$.  Let\/ $f\in\cG$
   be a reversible CR element, with a reversor\/ $r$ of order\/ $4$.
   Then, $r$ is conjugate to the linear map $($from\/ $\cA\setminus\cB
   )$ defined by the matrix\/ $R=\left(\begin{smallmatrix} 0 & -1 \\ 1
   & 0 \end{smallmatrix}\right)$, hence\/ $r^2$ is conjugate to\/
   $R^2 = I$.
\end{theorem}
\begin{proof}
Since $r^2$ is an involution and commutes with the CR element $f$,
Theorem~\ref{symnormform} shows that $r^2$ must be conjugate to $I =
\diag(-1,-1)$. Since $r$ itself is of finite order, it must be
conjugate to an element of $\cA$ or $\cE$, by
Fact~\ref{conjugate}. However, using the formulae given in the proof
of Lemma~\ref{fin-ord} for $e^n$ when $e\in\cE$, one deduces that
there can be no genuine order $4$ element that is elementary if
primitive $4$-th roots of unity are absent. In particular, this
excludes $e\in\cB$, see also the Appendix. So, $r$ is conjugate 
to an element $(\bs{a},M)\in \cA\setminus\cB$, with $M^2=I$.

Clearly, the matrix $R$ from the statement satisfies $R^2=I$, so it
is a root of $I$ in ${\rm GL}(2,K)$. Moreover, all other roots of $I$
in ${\rm GL}(2,K)$ are conjugate to $R$ in ${\rm GL}(2,K)$. To see
this, observe first that any $M\in {\rm GL}(2,K)$ with $M^2 = I$ must
satisfy $\trace(M)=0$ and $\det(M)=1$. This follows from a simple
direct calculation, which uses that $x^2=-1$ has no solution in $K$.
So, all solutions share the characteristic polynomial $P(x)=\det(x-M)
= x^2 + 1$.  This polynomial is irreducible over $K$ (by the
assumption on $U$), but splits as $P(x)=(x-i)(x+i)$ over the algebraic
closure $\hat{K}$ of $K$, with $i$ being a root of $-1$ in $\hat{K}$,
which cannot be in $U$ and hence not in $K$. This implies that $P(x)$
is also the minimal polynomial of all the possible
solutions. Consequently, they all have the same polynomial invariants,
hence are similar to one another, and also to the Frobenius companion
matrix of $P(x)$, which is the matrix $R$ (see \cite[Ch.~4.4]{AW} for
details).

Returning now to $r$, we have $r = h (\bs{a},M) h^{-1}$ for some
$h\in\cG$, $\bs{a}\in K^2$ and $M\in {\rm GL}(2,K)$ with $M^2=I$.
Since $1$ is not in the spectrum of $M$, $1-M$ is invertible. With
$\bs{c} = (1-M)^{-1}\bs{a}$, it is easy to check that
\[
   (\bs{c},1) (0,M) (-\bs{c},1) \; = \; (\bs{a},M)
\]
so that $(\bs{a},M)$ is conjugate, in $\cA$, to the linear map defined
by $M$. Now, putting things together, $r$ is conjugate to $M$ within
$\cG$ and, possibly employing one more ${\rm GL}(2,K)$-conjugation,
also to the linear map defined by the matrix $R$, as claimed.
\end{proof}

{}From Theorem~\ref{rev-order-real}, we can derive the possible structures
of $\cR(f)$, e.g., for $K=\RR$.

\begin{coro} \label{revgroup}
   Let\/ $K$ be a field of characteristic\/ $0$, with\/
   $U_K\simeq C_2$. If\/ $f\in\cG$ is a reversible CR element,
   $\cR(f)$ is one of the groups\/ $D_{\infty}\simeq C_{\infty}\rtimes C_2$,
   $C_{\infty}\rtimes C_4$, or\/ $(C_{\infty}\times C_2)\rtimes C_2$
   $($the last group comprising two different cases\/$)$.
\end{coro}
\begin{proof}
By Theorem~\ref{involutions} and Corollary~\ref{realisAbelian}, we
have either $\cS(f)\simeq C_{\infty}$ or $\cS(f)\simeq C_2\times
C_{\infty}$.  If $\cS(f)\simeq C_{\infty}$, a reversor $r$ of $f$ must
be an involution, whence $\cR(f)\simeq C_{\infty}\rtimes C_2 \simeq
D_{\infty}$.

Let $\cS(f)\simeq C_2\times C_{\infty}$ with involutory symmetry $s$,
which is then unique by Theorem~\ref{involutions}, and $C_\infty =
\langle h \rangle$. If the reversor $r$ is an involution, one has
$\cR(f)\simeq \cS(f)\rtimes C_2$.  Since $rsr$ is also an involutory
symmetry, we get $rsr=s$ by uniqueness, and $r$ and $s$ commute. Since
$r\neq s$, this gives $\cR(f) \simeq (C_{\infty}\times C_2)\rtimes
C_2$, with either $rhr^{-1}=h^{-1}$ (then giving $\cR(f)\simeq
D_\infty \times C_2$) or $rhr^{-1} = s h^{-1}$ (in which case $f$
must be an even power of $h$). Note that, in the latter case,
$\varrho = hr$ is an element of order $4$, and a reversor for $f$.

If $f$ has a reversor $r$ of order $4$, $r^2$ is an involutory
symmetry of $f$, hence unique and conjugate to $I=R^2$ with $R$ of
Theorem~\ref{natureorder4}.  This implies $\cS(f)\simeq C_2\times
C_{\infty}$ by Corollary~\ref{realisAbelian}, with $C_2 = \langle r^2
\rangle$, $C_{\infty} = \langle h \rangle$ and $f=r^{2\epsilon} h^m$
for $\epsilon\in\{0,1\}$ and some integer $m\ne 0$. In particular,
$r^2$ and $h$ commute, and $rhr^{-1}$ is a symmetry of $f$, so that
$rhr^{-1} = r^{2k} h^{\ell}$ for $k\in\{0,1\}$ and some
$\ell\in\ZZ$. Clearly, in view of $rfr^{-1} = f^{-1}$, this forces
$\ell=-1$.

If $k=0$, $r$ is also a reversor for $h$, and we have $\cR(f)\simeq
C_{\infty}\rtimes C_4$. This is the only case for $m$ odd, while for
$m$ even also $k=1$ is possible, i.e., $rhr^{-1} = r^{2} h^{-1}$. This
gives a group with the presentation 
\[
   \cR(f) \; = \; \langle r,h \mid r^4 = 1 , \; 
                  r h^{\pm 1} = h^{\mp 1} r^{-1} \rangle
\] 
which is an index $2$ extension of $\cS(f)\simeq C_{\infty}\times C_2$, 
but does not look like a simple semi-direct product. However, 
$\eta = h^{-1} r$ is an involution that satisfies
$\eta h \eta = r^2 h^{-1}$ and is a reversor for $f$.
This brings us back to $\cR(f)\simeq (C_{\infty}\times C_2)\rtimes
C_2$, where the outer $C_2$ is generated by $\eta$.
\end{proof}

Examples of CR elements $f\in\cG^{}_{\RR}$ of the generalised standard
form \eqref{defgenstan} illustrating all except the second possibility
of Corollary~\ref{revgroup} are given in \cite[Table 6]{RBstandard}
(in particular, one can extract examples for both subcases of
the third group).
To find an example of the remaining group structure (i.e.,
$C_{\infty}\rtimes C_4$), the simplest way \cite{UH} is to consider $f
= r e r e^{-1}$ with $e \! : \, x'=x+y^3, \, y'=y$ (which
commutes with $I$) and the linear map $r$ defined by the matrix $R$ of
Theorem~\ref{natureorder4}. Then, $f$ is reversible with reversor $r$,
but has no root in $\cG_{\RR}$ (though it has a root in $\cG_{\CC}$,
which then changes $\cS(f)$ and $\cR(f)$ in $\cG_{\CC}$).
The structure of this example will become more transparent from
Fact~\ref{factors} and Proposition~\ref{4normform} below.

\medskip
We now give various characterisations of reversible CR elements.  One
algebraic condition can be formulated via the poly-degree introduced in
\eqref{pol-deg-def}.  If we define the {\em reversal\/} of a finite
sequence of integers as $\overline{(n^{}_1,\dots,n^{}_k)} :=
(n^{}_k,\dots,n^{}_1)$, we observe
\begin{lemma} \label{degree-reversal}
   For all\/ $g\in\cG\setminus\cA$, one has\/
   $\pd(g^{-1}) = \overline{\pd(g)}$.
\end{lemma}
\begin{proof}
Assume $g$ is written in normal form. Its inverse is then a word in
affine and elementary mappings, potentially with an element from $\cB$
at the rightmost position. This gives a new sequence of degrees,
noting only those of the elementary mappings. Since $e$ and $e^{-1}$
have the same degree (compare \eqref{elementary} and
\eqref{inv-elementary}), for all $e\in\cE$, this new sequence is
nothing but $\overline{\pd(g)}$.

This sequence of degrees is not changed if the representation of $g^{-1}$ 
is now brought to normal form, by pulling the $\cB$-element to the left
and replacing, position by position, the mappings by the proper
representatives from $\cI$ and $\cJ$. So, $\overline{\pd(g)}$
is actually the poly-degree of $g^{-1}$, which proves the claim.   
\end{proof}

This enables us to formulate a rather restrictive necessary condition
for the reversibility of CR elements in $\cG$.
\begin{prop} \label{poly-reversing}
   Let the normal form of\/ $g\in\cG$ be cyclically reduced, which is then
   also true of the element\/ $g^{-1}$. A necessary condition for the
   reversibility of\/ $g$ is that\/ $\pd(g^{-1})$, which is the reversal
   of\/ $\pd(g)$, is a cyclic permutation of\/ $\pd(g)$.

   If, more generally, $g'$ is a CR element, it is conjugate to some
   element\/ $g$ with CRNF. The necessary
   condition for\/ $g'$ is then that the previous condition is met
   by\/ $g$. The outcome does not depend on the choice of\/ $g$.
\end{prop}
\begin{proof}
Let $g$ have a CRNF, which is then of length
$2m$ with $m\ge 1$, so that $\pd(g)$ is a sequence of length $m$
(recall that the poly-degree only keeps track of the elementary
maps). From Lemma~\ref{degree-reversal}, we know that $\pd(g^{-1}) =
\overline{\pd(g)}$, and the statement about $g$ now follows from the
result about conjugacy, see part~(3) of Proposition~\ref{conjugacy}.

If $g'$ is a CR element, we can't apply the criterion directly, but we
can pick any representative $g$ of the conjugacy class of $g'$ with
CRNF. Since $g'$ is reversible if and only
if $g$ is, the necessity of the claimed condition is obvious. It does
not depend on the choice of the representative because the poly-degrees
of different representatives are cyclic permutations of one another.
\end{proof}

\noindent {\sc Example}: 
Suppose $g'$ is a CR element, conjugate to a
$g$ in CRNF.  If $g$ contains up to two elements $e_i \in\cJ$, then
Proposition~\ref{poly-reversing} does not restrict $\pd(g)$ for $g$
(and $g'$) to be reversible (because any sequence of up to two
integers is a cyclic permutation of its reversal).  However,
restrictions generically arise when $g$ contains three or more
elements of $\cJ$. For instance, if $g$ has poly-degree $(2,3,4)$, it
can never be reversible.  This corresponds, in fact, to the lowest
degree of $g$ (i.e., $2\cdot 3\cdot 4=24$) for which $\pd(g)$ alone
can be exploited to rule out reversibility.

\medskip
We now proceed to describe, in more detail, the nature of reversible
elements of $\cG$, which will lead ultimately to the normal forms of
Proposition~\ref{invnormform} and Proposition~\ref{4normform} for
elements with involutory and order $4$ reversors, respectively.

\begin{prop} \label{det-condition}
   If\/ $g\in\cG$ has a reversor\/ $r\in\cG$, then\/ 
   $\det({\rm d}g) = \pm 1$.
\end{prop}
\begin{proof} 
By assumption, $g^{-1} = r g r^{-1}$ with $g,r\in\cG$.  Since the
Jacobians of polynomial automorphisms have constant determinant, a
simple application of the chain rule gives $\det({\rm d}g^{-1}) =
\det({\rm d}g)$, hence $\det({\rm d}g)^2 = 1$, which gives the claim.
\end{proof}

In view of the Remark after Proposition~\ref{normal-form},
reversibility puts an immediate restriction on the normal form.
\begin{coro}
   A necessary condition for\/ $g\in\cG$ to be reversible is that
   the element\/ $b\in\cB$ of its normal form\/ $(\ref{product-form})$ 
   satisfies\/ $\det({\rm d}b) = \pm 1$.   \qed
\end{coro}

Some further restrictions emerge for mappings which possess
fixed points.
\begin{prop}  \label{fixed}
   Let\/ $g$ be reversible, with reversor\/ $r$. If\/ $\bs{a}$ is a
   fixed point of\/ $g$, the Jacobian matrices\/ ${\rm d}g(\bs{a})$
   and\/ ${\rm d}g(r\bs{a})$ must have reciprocal spectrum.
\end{prop}
\begin{proof}
Since $g^{-1} = r\, g\, r^{-1}$, the chain rule (evaluated at the
point $r\bs{a}$) gives
\begin{equation*}
   {\rm d}g^{-1}(r \bs{a}) \; = \;
   {\rm d}r (\bs{a}) \, {\rm d}g (\bs{a}) \,
   {\rm d}r^{-1}(r \bs{a})\, .
\end{equation*}
Since $\,{\rm d}r (\bs{a})$ and ${\rm d}r^{-1}(r \bs{a})$ are the
inverses of each other (visible from the chain rule applied to $r^{-1}
r = 1$), ${\rm d}g^{-1}(r \bs{a})$ and ${\rm d}g (\bs{a})$ are
isospectral.

Observing $g^{-1} r\,\bs{a} = r\,g\,\bs{a} = r\,\bs{a}$ and applying
the chain rule to $g\,g^{-1} = 1$, one sees that ${\rm d}g^{-1}(r
\bs{a})$ is the inverse of ${\rm d}g (r \bs{a})$, from which the claim
follows.  
\end{proof}

To continue, we recall the following helpful factorisation property
from \cite{lamb}, formulated within the automorphism group of some
space. It will also shed more light on the examples discussed after
Corollary~\ref{revgroup}.
\begin{fact} \label{factors}
   An automorphism\/ $L$ is reversible, with reversor\/ $W\!$, if and
   only if some automorphism\/ $V$ exists such that\/ $L=V W^{-1}$
   together with\/ $V^2 = W^2$. In this case, also $V$ is a reversor.
\end{fact}
\begin{proof}
If $W$ is a reversor of $L$, define $V=LW$, which is invertible.
Clearly, $L=V W^{-1}$, and $V^2 = (LW)^2 = W^2$, as a consequence of
the relation $WLW^{-1}=L^{-1}$. Also, one quickly checks that
$VLV^{-1}=L^{-1}$.  Conversely, assuming $L=V W^{-1}$ with $V^2 =
W^2$, the last two relations follow immediately.  
\end{proof}

We now consider a normal form for reversible elements of $\cG$ which
have a reversing symmetry that is an involution.  Via
Fact~\ref{factors}, it follows that an automorphism is reversible with
an involutory reversor if and only if it is the product (i.e.,
composition) of two involutions (actually, this property goes back to
Birkhoff \cite{birk} whilst Fact~\ref{factors} represents a
generalisation of it).  Specialising to automorphisms in $\cG$, recall
that we know from Proposition~\ref{elementary-to-linear} that
involutions are conjugate to one of two possibilities:
$I=\diag(-1,-1)$ or $S=\diag(-1,1)$, equivalently
$T=\left(\begin{smallmatrix} 0 & 1 \\ 1 & 0
\end{smallmatrix}\right)$.
There are advantages to taking $T$ over $S$ in normal forms since the
former is in $\cA\setminus\cB$, indeed is in $\cI$ of \eqref{I-def}.
Note that $I$ is orientation-preserving, whereas $S$ and $T$ are
orientation-reversing. The canonical case of reversibility is that of
area-preserving maps which are the composition of two
orientation-reversing involutions. But the following result covers all
possibilities, not just this one.

\begin{prop} \label{invnormform}
   Let\/ $K$ be a field with $\ch(K)\neq 2$.  A CR element\/ $f\in\cG$ is
   reversible with a reversor that is an involution if and only if\/ $f$
   is conjugate to one of the following types of cyclically reduced
   normal forms:
\begin{equation} \label{invform1}
  \tilde{e}^{}_m \circ a^{}_{m-1}
  \circ {e}^{}_{m-1} \circ\ldots\circ a^{}_1
  \circ {e}^{}_1 \circ T \circ {e}_1^{-1} \circ a^{-1}_1
  \circ\ldots\circ {e}^{-1}_{m-1} \circ a^{-1}_{m-1} 
  \circ \tilde{e}^{-1}_m
\circ T
\end{equation}
\begin{equation} \label{invform2}
  \tilde{e}^{}_m \circ a^{}_{m-1}
  \circ {e}^{}_{m-1} \circ\ldots\circ a^{}_1
  \circ \hat{e} \circ a^{-1}_1
  \circ\ldots\circ {e}^{-1}_{m-1} \circ a^{-1}_{m-1} 
  \circ \tilde{e}^{-1}_m \circ T
\end{equation}
\begin{equation} \label{invform3}
  {a}^{}_m \circ e^{}_{m-1}
  \circ {a}^{}_{m-1} \circ\ldots\circ {e}^{}_1
  \circ {a}^{}_1 \circ \hat{e} \circ {a}_1^{-1} \circ e^{-1}_1
  \circ\ldots\circ {a}^{-1}_{m-1} \circ {e}^{-1}_{m-1} \circ
  a^{-1}_m \circ \bar{e}
\end{equation}
   In these normal forms, $T = \left(\begin{smallmatrix} 0 & 1 \\ 1 & 0
   \end{smallmatrix}\right)$, $a_i \in \cI$ of\/ $\eqref{I-def}$,
   $e_i \in \cJ$ of\/ $\eqref{J-def}$, $\hat{e}$ and $\bar{e}$ are
   particular cases of the involutions in $\cE\setminus\cB$ of the
   form \eqref{case-1} or \eqref{case-2}, and $\tilde{e}_m=b \circ
   e_m$, with $e_m \in \cJ$ and $b \in \cB$ a special case of
   \eqref{basic} $($as described further below\/$)$. Each normal form has
   $\det=\pm 1$ depending on the involutions present.  In each normal
   form, the only restriction on the appearance of $a_i$'s and $e_i$'s 
   is that an $e_i$ must occur if there is no elementary involution
   present, plus the form must be cyclically reduced.  It follows that
   if\/ $f'$ is any cyclically reduced element conjugate
    to\/ $f$, then\/ $\pd(f')= (n^{}_m,\ldots,n^{}_1,\hat{n}, n^{}_1,
    \ldots,n^{}_m, \bar{n})$, with all entries $\ge 2$ and $\hat{n}$
    and $\bar{n}$ absent or present according to the type of normal
    form above.
\end{prop}
\begin{proof}
{}From the fact that $f$ can be written as a composition of
involutions, together with Proposition~\ref{elementary-to-linear}, we
have that $f \in \cG$ is reversible with an involutory reversor if and
only if $f=h_1 T_1 h_1^{-1} h_2 T_2 h_2^{-1}$ with $h_i \in \cG$ and
$T_1, T_2 \in \{T,I\}$. Hence, the conjugate of $f$ given by $h_2^{-1}
f h_2$ takes the form
\begin{equation} \label{2form}
        h T_1 h^{-1} T_2,
\end{equation}
with $h=h_2^{-1} h_1$. We can take $h$ to be in the form
\eqref{product-form}. The cyclically reduced normal forms follow from
working through the possible forms of $h$, and combinations of $T_1$
and $T_2$. When $T_1$ and $T_2$ are different, it suffices to consider
$T_1=I$ and $T_2=T$, since the reverse possibility is conjugate to it.
A guiding principle, since $f$ is assumed to be CR, is that there
always remains an element of $\cE\setminus\cB$ and an element of
$\cA\setminus\cB$ in $h T_1 h^{-1} T_2$ after any possible reductions
into its cyclically reduced form.  Furthermore, we need the following
characterisations of involutions which follow from the proofs of
Lemma~\ref{affine-involutions} and
Proposition~\ref{elementary-to-linear} 
(with $\tilde{b} \in \cB$ in (i)-(iii)):
                                                                         
(i) an involution $a \in \cA\setminus\cB$ can be written
$a=\tilde{b}^{-1} \,T \,\tilde{b}$ with $\tilde{b}$ of the form
$x'=\alpha x + \gamma y + u$ and $y'=y+u$.
                                                                               
(ii) an involution $e \in \cE\setminus\cB$ conjugate to $T$ or $S$ can
be written as $e=\tilde{b}^{-1} \,\tilde{e} \,\tilde{b}$ with
$\tilde{e} \in \cE\setminus\cB$ of the form $x'=-x + y^2 \, Q(y), \,
y'=y$ or $x'=x + P(y), \, y'=-y$ with $Q \neq 0$ and $P(y)$ odd of
degree $\ge 3$.
                                                                               
(iii) an involution $e \in \cE\setminus\cB$ conjugate to
$I=\diag(-1,-1)$ can be written as $e=\tilde{b}^{-1} \tilde{e}
\tilde{b}$ with $\tilde{e} \in \cE\setminus\cB$ of the form 
$x'=-x + P(y), \, y'=-y$ with $P(y)$ even of degree $\ge 2$.
                                                                               
(iv) an involution in $\cB$ is conjugate in $\cB$ to one of $S$, $-S$
or $I$.
                                                                               
\noindent Note that in (ii)-(iv), $\tilde{b}$ of \eqref{basic} has
$\alpha=\beta=1$.
      
\smallskip 
We illustrate the reduction first for $T_1=T_2=T \in \cI$.  Firstly,
suppose $h$ ends in $e_1$ and begins with $e_m$, $m \ge 1$.  Then,
\eqref{2form} is the cyclically reduced element \eqref{invform1} with
$\tilde{e}_m$ actually in $\cJ$. If $h$ ends in $e_1$ and begins with
$b \circ a_m$, $a_m$ possibly missing, we conjugate \eqref{2form} and
consider $g=a_m^{-1} b^{-1} \,(h T h^{-1} T)\, b a_m$. This word ends
with the affine involution $a=a_m^{-1} b^{-1}\, T \,b a_m$. If this
involution belongs to $\cA\setminus\cB$, use characterisation (i)
above to see that $\tilde{b} g \tilde{b}^{-1}$ takes the form
\eqref{invform1}, with $\tilde{e}_m=\tilde{b} e_m$.  Now $\tilde{e}_m$
is in $\cE\setminus\cB$, but possibly not in $\cJ$. Otherwise, the
involution is $a=a_m^{-1} b^{-1} T b a_m \in \cB$, conjugate in $\cB$
to $\pm S$ from (iv) above. Then, consider $g'=e^{-1}_m g e_m$ which
ends with the elementary involution $e=e^{-1}_m a e_m$. If $e \in
\cE\setminus\cB$, use characterisation (ii) above to see that $g'$ is
conjugate to \eqref{invform2} with $\hat{e}$ of one of the forms
described.  Whereas, if $e \in \cB$, one continues by considering
$g''=a_{m-1}^{-1} g' a_{m-1}$, which ends in the affine involution
$a'=a_{m-1}^{-1} e a_{m-1}$ which is either in $\cA\setminus\cB$ or in
$\cB$. It is clear how this repeated process must eventually exhaust
itself.
                                                                               
Next, suppose $h$ in \eqref{2form} takes the form $h= e^{}_m
\circ\ldots\circ e^{}_1 \circ a^{}_1$, i.e., $h$ ends with an affine
coset representative $a_1$ and begins with $e_m$, $m \ge 1$. Then,
\eqref{2form} contains the affine involution $a=a_1 T a^{-1}_1$. If $a
\in \cA\setminus\cB$, use characterisation (i) again and rewrite $h$
in the form \eqref{product-form} to obtain
\eqref{invform1}. Otherwise, $a \in \cB$ is conjugate to $\pm S$ and
one moves on to study the elementary involution $e=e^{}_1 a e^{-1}_1$.
This process leads to a cyclically reduced word (\ref{invform2}) with
$\hat{e}$ of characterisation (ii) above if $e \in \cE\setminus\cB$,
or returns once more to the study of an affine involution $a_2 e
a^{-1}_2$ etc.  

{}Finally, consider the case that $h$ ends in $a_1 \in \cI$ but begins
with $b \circ a_m$, whence we have $h= b \circ a^{}_m \circ e^{}_{m-1}
\circ \ldots \circ e^{}_1 \circ a^{}_1$.  Now one uses, in tandem, the
combination of the above-mentioned procedures. One takes $g=a_m^{-1}
b^{-1} h T h^{-1} T b a_m$ and sees that the processes will exhaust
themselves in one of (\ref{invform1})--(\ref{invform3}), with the
elementary involutions occurring being those of characterisation (ii).
                                                                               
The cases in \eqref{2form} when $T_1=T_2=I$ and when $T_1=I$ and
$T_2=T$ follow a similar, but simpler, path.  This is because linear
elements such as $a_i$ commute with $I$.  This leads to less cases
that need to be considered. When $T_1=T_2=I$, we obtain
\eqref{invform3} with both $\hat{e}$ and $\bar{e}$ of the form
described in characterisation (iii) above. When $T_1=I$ and $T_2=T$,
we obtain \eqref{invform2} or \eqref{invform3} with $\hat{e}$ of
characterisation (iii) and $\bar{e}$ of characterisation (ii) above.
                                                                               
We remark that, without loss of generality, the element $b \in \cB$
occurring at the start of (\ref{invform1})--(\ref{invform3}) can be
chosen from the quotient of $\cB$ and the centraliser in $\cB$ of the
last element of (\ref{invform1})--(\ref{invform3}). For example, for
(\ref{invform1})--(\ref{invform2}), this gives an element
\eqref{basic} containing just $3$ parameters instead of $5$.
\end{proof}
  
\noindent {\sc Remark}: The normal forms of
Proposition~\ref{invnormform} are similar to those found in
\cite[Thm.~1]{GM}.  There, the authors express their
cyclically reduced normal forms using compositions of H\'enon maps $h_m
\circ \ldots \circ h_1$ (and the
inverse of such a composition) with $h_i$
of $\eqref{defgenhenonmap}$, instead of our
expressions above in terms of $a_i$ and $e_i$.

\medskip
Finally, we consider a normal form for reversible
elements of $\cG$ with a reversor of order\/ $4$.

\begin{prop} \label{4normform}
   Let $K$ be a field with $\ch(K)\neq 2$, with a unit group\/ $U$
   that contains\/ $\{\pm 1\}$, but no primitive $4$-th root of unity
   $($thus including the case\/ $U\simeq C_2$\/$)$.  A CR element
   $f\in\cG$ is reversible with a reversor $r$ of order\/ $4$ if and
   only if f is conjugate to the CRNF
\begin{equation} \label{4form1}
      {e}^{}_m \circ\ldots\circ a^{}_1
      \circ e^{}_1 \circ R^{}_1 \circ e_1^{-1} \circ a^{-1}_1
      \circ\ldots\circ e^{-1}_m \circ R^{}_2 \ts .
\end{equation}
   Here, $R_1 = R = \left(\begin{smallmatrix} 0 & -1 \\ 1 & 0
   \end{smallmatrix}\right) \in   \cA\setminus\cB$
   and $R_2= \left(\begin{smallmatrix} \alpha &
   -(\alpha^2+1)/\gamma \\ \gamma & -\alpha \end{smallmatrix}\right)
   \in \cA\setminus\cB$ $($since $\gamma \neq 0)$, the latter including
   $R_2 = R$ via $\alpha=0$ and $\gamma=1$. Moreover, $m \ge
   1$, $e_i \in \cJ$ of\/ \eqref{J-def} must have\/ $P_i(y)$ odd and
   $e^{}_1$ and $e^{}_m$ must appear.  It follows that $f$ necessarily
   has\/ $\det({\rm d}f) = 1$ and has a fixed point.
    Also, if $f'$ is a cyclically reduced element conjugate to $f$,
    $\pd(f')= (n^{}_m,\ldots,n^{}_1, n^{}_1,\ldots,n^{}_m)$ where\/
    $n_{i}$ are odd integers $\ge 3$, and $f'$ commutes with
    $x'=-x+u$, $y'=-y+v$, with some $u,v\in K$.
\end{prop}
\begin{proof}
{}From Fact~\ref{factors}, one can see that an automorphism $L$ has a
reversor $W$ of order $4$ if and only if $L = V W^{-1}$ with $V^2=W^2$
and $V$ also a reversor of order $4$.  Take in \eqref{4form1}
$W^{-1}=R_2$ and for $V$ the first term, conjugate to $R_1$. We see
that $V^2=W^2=\diag(-1,-1)$ for either possibility of $R_1$ and $R_2$
under the assumptions given on $e_m$.  Hence, \eqref{4form1} has order
$4$ reversors, e.g., $V$ and $W=R_2^3=-R_2$, and this
property will be preserved under conjugacy.

Consider the converse. Since Theorem~\ref{natureorder4} characterises
the order $4$ reversors, we have, using Fact~\ref{factors}, that
$f=h_1 R h_1^{-1} h_2 R h_2^{-1}$ with $h_i\in \cG$ and $R$ the matrix
in the statement (note that $W$ an order $4$ reversor implies the same
property for $V$ and $W^{-1}$).  This shows immediately that
$\det({\rm d}f)=1$.  Hence, the conjugate of $f$ given by $h_2^{-1} f
h_2$ takes the form
\begin{equation} \label{4form}
        h R h^{-1} R
\end{equation}
with $h=h_2^{-1} h_1$.
                                                                               
It follows from Fact~\ref{factors} that $h$ (and hence \eqref{4form})
commutes with $R^2=I=\diag(-1,-1)$, equivalently $I\ts h\ts I = h$. If
we take for $h$ an expression of the form \eqref{product-form},
identical reasoning to that used in the proof of
Theorem~\ref{symnormform} establishes that $h$ has $b$ linear with
$e_i$ having $P_i(y)$ odd.  In the expression for $h$, there must be
at least one $e_i$, otherwise \eqref{4form} $\in \cA$ and $f$ is then
conjugate to an affine element, in contradiction to $f$ being CR.
                                                                               
Next, we need to consider the different possibilities for $h$ and the
reduction of \eqref{4form}, if necessary, to a cyclically reduced
word. If $h$ ends with an element from $\cJ$, we obtain \eqref{4form1}
with $R_1=R$.  If, in addition, $h$ starts with $b \circ a_m$ followed
by an element from $\cJ$, $a_m$ possibly absent and $b$ possibly the
identity, study the conjugate $a_m^{-1} b^{-1} \,(h R h^{-1} R)\, b
a_m$.  It is of the form \eqref{4form1} and ends with a linear
traceless order 4 element $R_2=a_m^{-1} b^{-1} R b a_m$ of the form
indicated. A straightforward calculation shows that the entry $\gamma$
is necessarily non-zero, because $x^2=-1$ has no solution in $K$ by
assumption.

Otherwise, if the last element of $h$ was $a \in \cI$, we could write
$h R h^{-1} = h' (a R a^{-1}) {h'}^{-1}$, with $h'$ ending in an
elementary map and with the new linear order $4$ element $R'=a R
a^{-1}= \left(\begin{smallmatrix} -\beta & 1 \\ -(1+\beta^2) & \beta
\end{smallmatrix}\right)$.  Again, $1+\beta^2 \neq 0$ in $K$ by
assumption, so $R' \in \cA\setminus\cB$. However, we can rewrite $R'=b
R b^{-1}$ with $b= \left(\begin{smallmatrix} 1 & -\beta \\ 0 &
-(1+\beta^2) \end{smallmatrix}\right) \in \cB$. Hence $h R
h^{-1}=(h'b) R (h'b)^{-1}$, where $h' b$ takes the form
\eqref{product-form} ending in an element from $\cJ$. This returns us
to the case of the previous paragraph.  Considering now the start of
$h' b$ as above, and possibly using a further affine conjugacy, again
returns the form \eqref{4form1} with $R_1=R$ and $R_2$ as given. 

This explains the normal form given. The result for $\pd(f')$ is a
direct consequence of the odd nature of the $P_i(y)$ in $e_i$. The
fact that $f'$ has the symmetry indicated follows from the
fact that $R_i^2=I=\diag(-1,-1)$ commutes with \eqref{4form1} and from
Corollary~\ref{involtest}.
\end{proof}


\noindent {\sc Remark}:
For the case $K=\CC$, reference \cite{GM2} presents normal forms for
CR elements that possess a reversor of order $2n$.

\section*{Appendix: Elements of $\cB$ of finite order}

Symmetries of CR elements of finite order are conjugate to elements of
$\cB$ of finite order. As these are of particular relevance for
detecting existing symmetries, we add a short classification here, for
an arbitrary field $K$.

Recall that $\cB = \{ (\bs{a}, M) \mid \bs{a}\in K^2, \, M \in \cT\} =
K^2 \rtimes \cT$ where $\cT$ denoted the subgroup of all upper
triangular matrices of ${\rm GL}(2,K)$. Since
\begin{equation*}
   (\bs{a}, M)^n \; = \; \big( (1+M+M^2+\ldots +M^{n-1})\bs{a}, M^n\big),
\end{equation*}
it is clear that $(\bs{a}, M)^n = (\bs{0}, 1)$ implies $M^n=1$ and
$(1+M+M^2+\ldots +M^{n-1})\bs{a} = \bs{0}$.

Consider a matrix $M=\bigl( \begin{smallmatrix} \alpha & \gamma \\
 0 & \beta \end{smallmatrix} \bigr)$ with $\alpha,\beta,\gamma\in K$
and $\alpha\beta\neq 0$, so that $M$ is invertible.
\begin{lemma}
   For\/ $n\in\ZZ$, the matrix powers of\/ $M$ are
\begin{equation*}
   M^n \; = \; \begin{pmatrix} \alpha^n & \gamma(n) \\
               0 & \beta^n  \end{pmatrix}
\end{equation*}
   where\/ $\gamma(0)=0$, and, for all\/ $n\ge 1$,
   $\gamma(-n) = -\gamma(n)/(\alpha\beta)^n$ with
\begin{equation*}
   \gamma(n) \; = \; \gamma \sum_{m=0}^{n-1} \alpha^m \beta^{n-1-m} .
\end{equation*}
\end{lemma}
\begin{proof}
The formula for $\gamma(n)$, for positive $n$, is easy to check by
induction, while the inversion formula for $2\!\times\! 2$-matrices
gives the result for negative $n$, and $\gamma(0)=0$ is clear.  
\end{proof}

If $M^n = 1$, we must have $\alpha^n = \beta^n = 1$ and $\gamma(n)=0$.
If $\alpha^n = \beta^n = 1$, but $\alpha\neq\beta$, a simple geometric
series argument shows that $\gamma(n)=0$ is automatic. On the other
hand, if $\alpha=\beta$, one finds
\begin{equation*}
   \gamma(n) \; = \; n \, \alpha^{n-1} \gamma \, .
\end{equation*}
In characteristic $0$, this can only vanish for $\gamma=0$. Otherwise,
$\gamma(n)$ vanishes also if $\ch(K)$ divides $n$. Consequently,
$\ord(M) = \lcm(\ord(\alpha),\ch(K))$. This gives:
\begin{prop}  \label{finord-matrices}
   Consider\/ $M=\bigl( \begin{smallmatrix} \alpha & \gamma \\
   0 & \beta \end{smallmatrix} \bigr)$ with $\alpha,\beta\in U_K$,
   and let $n=\lcm(\ord(\alpha),\ord(\beta))$.
   If\/ $\ch(K)=0$, $M$ is of finite order if either\/ $\alpha\neq\beta$
   or\/ $\alpha=\beta$ with\/ $\gamma=0$. In both cases, $\ord(M)=n$.
   If\/ $\ch(K)\neq 0$, $M$ is of finite order for all\/ $\gamma\in K$, 
   with\/ $\ord(M)=n$ for\/ $\alpha\neq\beta$ and\/
   $\ord(M)=\lcm(\ch(K),n)$ for\/ $\alpha = \beta$.   \qed
\end{prop}

Now, we have to extend to the affine case. Let $M$ be a matrix of
order $n$. If $(1-M)$ is invertible, another geometric series argument
shows that all affine extensions $(\bs{a},M)$ are also of order $n$.
So, assume $(1-M)$ is not invertible, i.e., $1$ is an eigenvalue of
$M$.  If $M\bs{x} = \bs{x}$, with $\bs{x}\neq\bs{0}$, one has
$(1+M+\ldots +M^{m-1})\bs{x} = m \bs{x}$, which vanishes only for
$\ch(K) | m$. This always happens if $m$ is some multiple of $n$, as
long as $\ch(K)\neq 0$.  In characteristic $0$, however, the
translational part of the affine extension has to avoid the kernel of
$(M-1)^k$, where $k$ is the exponent of the factor $(x-1)$ in the
minimal polynomial of $M$.

\begin{prop}  \label{finord-affine}
   Let\/ $M\in\cT$ with\/ $\ord(M)=n < \infty$. Let\/ $k$ be the
   exponent of\/ $(x-1)$ in the minimal polynomial of\/ $M$, and set\/
   $S = \ker\bigl( (M-1)^k\bigr)$.  If\/ $\ch(K)=0$, the element\/
   $(\bs{a},M)$ with\/ $\bs{a}\in K^2$ is of finite order iff\/
   $\bs{a}$ has no component in the generalised eigenspace\/ $S$. In
   this case, the order is\/ $n$.  If\/ $\ch(K)\neq 0$, $(\bs{a},M)$
   is of finite order for all\/ $\bs{a}\in K^2$, but the order can be
   a multiple of\/ $n$.  \qed
\end{prop}

\section*{Acknowledgment}

It is a pleasure to thank Ulrich Hermisson for his cooperation and for
helpful discussions.  M.B.\ would like to thank the School of
Mathematics of the University of New South Wales for financial support
during a stay in November 2002, where a substantial part of this work
was done.


\bigskip

\begin{thebibliography}{99}
\small

\bibitem{AW}
W.\thinspace A.\ Adkins and S.\thinspace H.\ Weintraub,
\textit{Algebra -- An Approach via Module Theory},
corr.\ printing, Springer, New York (1999).

\bibitem{AR}
P.~Ahern and W.~Rudin,
\textit{Periodic automorphisms of\/ $\CC^n$},
Indiana Univ.\ Math.\ J.\ {\bf 44} (1995) 287--303.

\bibitem{Asa}
T.\ Asanuma,
\textit{Non-linearizable algebraic group action on $\bs{A}^n$},
J.\ Algebra {\bf 166} (1994) 72--79.

\bibitem{birk}
G.\thinspace D.\ Birkhoff,
\textit{The restricted problem of three bodies},
Rend.\ Circ.\ Mat.\ Palermo {\bf 39} (1915) 265--334.

\bibitem{BRcat}
M.\ Baake and J.\thinspace A.\thinspace G.\ Roberts,
\textit{Reversing symmetry group of $GL(2,\ZZ)$ and $PGL(2,\ZZ)$
matrices with connections to cat maps and trace maps},
J.\ Phys.\ A:\ Math.\ Gen.\ {\bf 30} (1997) 1549--1573.

\bibitem{BRtorus}
M.\ Baake and J.\thinspace A.\thinspace G.\ Roberts,
\textit{Symmetries and reversing symmetries of toral automorphisms},
Nonlinearity {\bf 14} (2001) R1--R24; math.DS/0006092.

\bibitem{Cohen}
D.\thinspace E.\ Cohen,
\textit{Combinatorial Group Theory:\ A Topological Approach},
Cambridge University Press, Cambridge (1989).

\bibitem{CM}
H.\thinspace S.\thinspace M.\ Coxeter and 
W.\thinspace O.\thinspace J.\ Moser,
\textit{Generators and Relations for Discrete Groups},
4-th ed., Springer, Berlin (1980).

\bibitem{Essen}
A.~van den Essen,
\textit{Seven lectures on polynomial automorphisms}, in:\
\textit{Automorphisms of Affine Spaces}, ed.\ A.\ van den Essen,
Kluwer, Dordrecht (1995), pp.\ 3--39.

\bibitem{Essen2}
A.~van den Essen,
\textit{Polynomial Automorphisms and the Jacobian Conjecture},
Birkh\"auser, Basel (2000).

\bibitem{FM}
S.~Friedland and J.~Milnor,
\textit{Dynamical properties of plane polynomial automorphisms},
Ergod.\ Th.\ \& Dynam.\ Syst.\ {\bf 9} (1989) 67--99.

\bibitem{GM}
A.~G\'omez and J.\thinspace D.~Meiss,
\textit{Reversible polynomial automorphisms of the plane:\
the involutory case}, 
Phys.\ Lett.\ {\bf A 312} (2003) 49--58; nlin.CD/0209055.

\bibitem{GM2}
A.~G\'omez and J.\thinspace D.~Meiss,
\textit{Reversors and symmetries for polynomial automorphisms
of the complex plane},
Nonlinearity {\bf 17} (2004) 975--1000;
nlin.CD/0304035.

\bibitem{Goodson}
G.\thinspace R.~Goodson,
\textit{Inverse conjugacies and reversing symmetry groups},
Amer.\ Math.\ Monthly {\bf 106} (1999) 19--26.

\bibitem{UH}
U.~Hermisson, private communication (2002).

\bibitem{Jung}
H.\thinspace W.\thinspace E.\ Jung,
\textit{\"Uber ganze birationale Transformationen der Ebene},
J.\ Reine Angew.\ Math.\ (Crelle) {\bf 184} (1942) 161--174.

\bibitem{Kamba}
T.\ Kambayashi,
\textit{Automorphism group of a polynomial ring and algebraic
group action on an affine space},
J.\ Algebra {\bf 60} (1979) 439--451.

\bibitem{KS}
A.~Karrass and D.~Solitar,
\textit{The subgroups of a free product of two groups with
an amalgamated subgroup},
Trans.\ AMS {\bf 150} (1970) 227--255.

\bibitem{Kraft}
H.\ Kraft and G.\ Schwarz,
\textit{Finite automorphisms of affine space}, in:\
\textit{Automorphisms of Affine Spaces}, ed.\ A.\ van den Essen,
Kluwer, Dordrecht (1995), pp.\ 55--66.

\bibitem{Kulk}
W.~van der Kulk,
\textit{On polynomial rings in two variables},
Nieuw Arch.\ Wisk.\ {\bf 1} (1953) 33--41.

\bibitem{lamb}
J.\thinspace S.\thinspace W.\ Lamb,
\textit{Reversing symmetries in dynamical systems},
J.\ Phys.\ A:\ Math.\ Gen.\ {\bf 25} (1992) 925--937.

\bibitem{LR}
J.\thinspace S.\thinspace W.\ Lamb and 
J.\thinspace A.\thinspace G.\ Roberts,
\textit{Time-reversal symmetry in dynamical systems:\ A survey},
Physica {\bf D 112} (1998) 1--39.

\bibitem{Lang}
S.~Lang,
\textit{Algebra}, rev.\ 3rd ed., Springer, New York (2002).

\bibitem{LS}
R.\thinspace C.~Lyndon and P.\thinspace E.~Schupp,
\textit{Combinatorial Group Theory},
Springer, Berlin (1977); reprint (2001).

\bibitem{MKS}
W.~Magnus, A.~Karrass and D.~Solitar,
\textit{Combinatorial Group Theory:\ Presentations of Groups
in Terms of Generators and Relations}, 2nd ed.,
Dover, New York (1976).

\bibitem{Mol}
D.\thinspace I.~Moldavanskii,
\textit{Certain subgroups of groups with one defining relation},
Sibirsk.\ Mat.\ $\check{\rm Z}$.\ {\bf 8} (1967) 1370--1384.

\bibitem{RBtrace}
J.\thinspace A.\thinspace G.~Roberts and M.~Baake, 
\textit{Trace maps as 3D reversible dynamical systems with an invariant},
J.\ Stat.\ Phys.\ {\bf 74} (1994) 829--888. 

\bibitem{RBstandard}
J.~A.~G.~Roberts and M.~Baake,
\textit{Symmetries and reversing symmetries of area-preserving 
polynomial mappings in generalised standard form},
Physica {\bf A 317} (2003) 95--112; math.DS/0206096.

\bibitem{roqu}
J.\thinspace A.\thinspace G.\ Roberts and 
G.\thinspace R.\thinspace W.\ Quispel,
\textit{Chaos and time-reversal symmetry --- order and chaos in reversible
dynamical systems},
Phys.\ Rep.\ {\bf 216} (1992) 63--177.

\bibitem{rovi1}
J.\thinspace A.\thinspace G.~Roberts and F.~Vivaldi,
\textit{Arithmetical method to detect integrability in maps},
Phys.\ Rev.\ Lett.\ {\bf 90} (2003) 034102.

\bibitem{rovi2}
J.\thinspace A.\thinspace G.~Roberts and F.~Vivaldi,
\textit{Signature of time-reversal symmetry in
polynomial automorphisms over finite fields}, 
preprint (2004).

\bibitem{rowi}
J.\thinspace A.\thinspace G.\ Roberts and 
R.\thinspace S.\ Wilson,
\textit{Reversibility of orientation-reversing cat maps 
and the amalgamated free product structure of\/ {\rm PGL(2,$\ZZ$)}},
preprint (2002).

\bibitem{Rudin}
W.~Rudin,
\textit{Injective polynomial maps are automorphisms},
Amer.\ Math.\ Monthly {\bf 102} (1995) 540--543.

\bibitem{Serre}
J.-P.~Serre,
\textit{Trees},
Springer, Berlin (1980); 2nd corr.\ printing (2003).

\bibitem{V1}
A.\thinspace P.~Veselov,
\textit{Integrable maps},
Russian Math.\ Surveys {\bf 46} (1991) 1--51. 

\bibitem{V2}
A.\thinspace P.~Veselov,
\textit{Growth and integrability in the dynamics of mappings},
Commun.\ Math.\ Phys.\ {\bf 145} (1992) 181--193.

\bibitem{Wright}
D.~Wright,
\textit{Abelian subgroups of\/ ${\rm Aut}_k (k[X,Y])$ and
applications to actions on the affine plane},
Illinois J.\ Math.\ {\bf 23} (1979) 579--634.

\bibitem{W2}
D.~Wright, 
private communication (2004).

\end{thebibliography}
\end{document}